\documentclass[reqno]{amsart}

\usepackage{amssymb,amscd,amsxtra,epic,eepic}

\numberwithin{equation}{section}
\newcommand\note[1]%
{$^\dagger$\marginpar{\sf\footnotesize{$^\dagger${#1}}}}

\swapnumbers
\newtheorem{theorem}{Theorem}[section]
\newtheorem{proposition}[theorem]{Proposition}
\newtheorem{lemma}[theorem]{Lemma}
\newtheorem{corollary}[theorem]{Corollary}
\newtheorem{addendum}[theorem]{Addendum}

\newtheorem{problem}[theorem]{Problem}

\theoremstyle{definition}
\newtheorem{definition}[theorem]{Definition}
\newtheorem{example}[theorem]{Example}
\newtheorem{remark}[theorem]{Remark}

\newcommand\eu{\mathfrak}
\newcommand\lie{\mathfrak}
\renewcommand\k{\lie{k}} 
\renewcommand\t{\lie{t}}
\newcommand\g{\lie{g}}
\newcommand\h{\lie{h}} 
\renewcommand\u{\lie{u}} 
\newcommand\s{\lie{s}} 
\newcommand\p{\lie{p}} 
\newcommand\q{\lie{q}} 
\renewcommand\a{\lie{a}}

\newcommand\bb{\mathbb}

\newcommand\Z{\bb{Z}} 
\newcommand\Q{\bb{Q}}
\newcommand\R{\bb{R}} 
\newcommand\C{\bb{C}}
\renewcommand\P{\bb{P}}

\newcommand\ca{\mathcal}

\newcommand\Ad{\operatorname{Ad}}
\newcommand\hull{\operatorname{hull}}

\renewcommand\star{\operatorname{star}}

\newcommand\grad{\operatorname{grad}}

\newcommand\Inn{\operatorname{Inn}}

\newcommand\Hom{\operatorname{Hom}}

\newcommand\id{\operatorname{id}}

\renewcommand\Re{\operatorname{Re}}
\renewcommand\Im{\operatorname{Im}}

\newcommand\Lie{\operatorname{Lie}} 

\newcommand\curv{\operatorname{curv}}

\newcommand\SO{\operatorname{\mathbf{SO}}}
\newcommand\GL{\operatorname{\mathbf{GL}}}

\newcommand\SU{\operatorname{\mathbf{SU}}}
\newcommand\U{\operatorname{\mathbf{U}}}
\renewcommand\Sp{\operatorname{\mathbf{Sp}}}

\newcommand\BC{\operatorname{\mathbf{BC}}}
\newcommand\CC{\operatorname{\mathbf{C}}}

\newcommand\suc{\succcurlyeq}

\newcommand\qu{/\kern-.7ex/}
\newcommand\bigqu{\big/\kern-.85ex\big/}

\newcommand\longhookrightarrow{\lhook\joinrel\longrightarrow}

\newcommand\map{\longrightarrow}
\newcommand\inj{\longhookrightarrow}
\newcommand\sur{\longrightarrow\kern-1.9ex\to}
\newcommand\iso{\longhookrightarrow\kern-1.9ex\to}

\newcommand\eps{\varepsilon}
\newcommand\bu{{\scriptscriptstyle\bullet}}
\newcommand\ti{\tilde}

\newcommand\antiddots{\mathinner{%
\mkern1mu\raise1pt\vbox{\kern7pt\hbox{.}}
\mkern2mu\raise4pt\hbox{.}\mkern2mu\raise7pt\hbox{.}\mkern1mu}}

\newcommand\inv{^{-1}} 

\renewcommand\subset{\subseteq}
\renewcommand\supset{\supseteq}

\newcommand\rel{_{\mathrm{rel}}}
\newcommand\com{_{\mathrm{com}}}
\newcommand\hor{_{\mathrm{hor}}} 
\newcommand\reg{_{\mathrm{reg}}}

\newcommand\dup{^{\mathrm{dup}}}

\title{Moment Maps and Riemannian Symmetric Pairs}

\author{Luis O'Shea}
\address{Department of Mathematics, Cornell University, Ithaca, New
York 14853-7901} 
\email{luis@math.cornell.edu}
\author[Reyer Sjamaar]{Reyer Sjamaar$^*$}
\thanks{$^*$Partially supported by an Alfred P. Sloan Research
Fellowship and by NSF Grant DMS-9703947}
\email{sjamaar@math.cornell.edu}
\date{8 February 1999.  Revised 25 October 1999}

\begin{document} 


\begin{abstract}
We study Hamiltonian actions of a compact Lie group on a symplectic
manifold in the presence of an involution on the group and an
antisymplectic involution on the manifold.  The fixed-point set of the
involution on the manifold is a Lagrangian submanifold.  We
investigate its image under the moment map and conclude that the
intersection with the Weyl chamber is an easily described subpolytope
of the Kirwan polytope.  Of special interest is the integral K\"ahler
case, where much stronger results hold.  In particular, we obtain
convexity theorems for closures of orbits of the noncompact dual group
(in the sense of the theory of symmetric pairs).  In the abelian case
these results were obtained earlier by Duistermaat.  We derive
explicit inequalities for polytopes associated with real flag
varieties.
\end{abstract}


\maketitle

\tableofcontents

\section{Introduction}\label{section;introduction}

The main result of this paper is a ``real form'' of Kirwan's convexity
theorem.  We apply our result to flag varieties of real semisimple
groups and thus obtain eigenvalue inequalities, which generalize
inequalities found by Weyl, Ky Fan, Kostant, Klyachko, and many
others.  For instance, we find an answer to the following question:
given the singular values of two rectangular matrices (of the same
size), what are the possible singular values of their sum?

We obtain these results by studying symplectic manifolds $M$ equipped
with an antisymplectic involution $\tau$ and a Hamiltonian action of a
compact Lie group $U$.  We assume the fixed-point set of the
involution to be nonempty, in which case it is a Lagrangian
submanifold.  We denote the moment map for the $U$-action on $M$ by
$\Phi$.  Our goal is to describe the image of the Lagrangian under
$\Phi$.  Good results can of course only be expected if the involution
and the action are related in a reasonable way.  The right condition
to impose is for $\Phi$ to be ``anti-equivariant'' relative to the
involution $\tau$ and an involutive automorphism $\sigma$ on the group
$U$.  The Lagrangian is then not stable under the action of $U$, but
only under the fixed-point group $U^\sigma$.  Thus the theory of
Riemannian symmetric pairs comes into play.

If $U$ is a torus and $\sigma$ is inversion, anti-equivariance amounts
to \emph{in}variance of $\Phi$ with respect to $\tau$.  This special
case was studied by Duistermaat \cite{duistermaat;convexity}, who
proved that the image of $M^\tau$ under $\Phi$ is equal to the full
image $\Phi(M)$.

This turns out to be false in general.  Our main result is the
following description of $\Phi(M^\tau)$ as a subset of $\Phi(M)$.
Instead of the image $\Phi(M)$, it is easier to consider the
intersection $\Delta(M)=\Phi(M)\cap\t^*_+$, where $\t^*_+$ is a Weyl
chamber in $\u^*$.  According to Kirwan, $\Delta(M)$ is a convex
polytope, called the \emph{moment polytope} of $M$.  Our result is
that the set $\Delta(M^\tau) =\Phi(M^\tau)\cap\t^*_+$ is equal to the
subpolytope of $\Delta(M)$ consisting of symmetric elements.  In other
words, $\Delta(M^\tau)$ is equal to the intersection of $\Delta(M)$
with the Weyl chamber of the \emph{restricted} root system.

The proof is an adaptation of the proof of Kirwan's convexity theorem
given in \cite{sjamaar;convexity}.  (We also attempted to adapt the
perhaps simpler approach followed by Lerman \emph{et al.}\
\cite{lerman;nonabelian-convexity}, but were unable to do so.)
Briefly, the method is to establish the result first in the K\"ahler
case (more precisely, for projective and affine varieties), and then
handle the general case by dint of a local normal form theorem, which
asserts that near the orbits that matter every Hamiltonian
$U$-manifold is isomorphic to an affine variety in a suitable sense.
An advantage of this approach is that it yields a characterization of
the vertices and an attractive, if unwieldy, invariant-theoretic local
description of the polytope near any point.  More importantly, the
K\"ahler case has the noteworthy feature that $G$, the noncompact real
reductive group dual to $U$, acts holomorphically and \emph{commutes}
with the antiholomorphic involution $\tau$.  We are thus led to a
convexity theorem for closures of $G$-orbits, which is a nonabelian
version of a result of Atiyah's.

In Section \ref{section;elementary} we introduce the category of
Hamiltonian $(U,\sigma)$-manifolds and give a number of examples.  Our
main result is stated in Section \ref{section;convex}.  Section
\ref{section;quantization} is an expos\'e of K\"ahler quantization for
Hamiltonian $(U,\sigma)$-manifolds, where we develop the techniques
needed for establishing the convexity theorem in the K\"ahler case,
which we carry out in Sections \ref{section;kahler} and
\ref{section;affine}.  We then state the normal form theorems in
Section \ref{section;normal} and finish the proof in Section
\ref{section;proof}.  Applications inspired by mechanics (cotangent
bundles) and representation theory (coadjoint orbits) are discussed in
Section \ref{section;application}.  A few elementary results on
antisymplectic maps and symmetric pairs are covered in Appendices
\ref{section;antisymplectic} and \ref{section;symmetric}.  Appendix
\ref{section;notation} is a list of notational conventions.

\section{Elementary properties and examples}\label{section;elementary}

Let $U$ be a compact Lie group.  Let $\sigma$ be an involution of $U$,
i.e.\ an automorphism which satisfies $\sigma^2=\id_U$.  (We are
mostly interested in the case where $\sigma$ is not the identity.)
The derivative of $\sigma$ induces involutions on the Lie algebra $\u$
and on the dual space $\u^*$, both of which we denote for simplicity
also by $\sigma$.

\begin{definition}
A \emph{Hamiltonian $(U,\sigma)$-manifold} is a quadruple
$(M,\omega,\tau,\Phi)$ consisting of a smooth manifold $M$ equipped
with an action of the group $U$, a symplectic form $\omega$, an
antisymplectic involution $\tau$ (i.e.\ a diffeomorphism $\tau$
satisfying $\tau^2=\id_M$ and $\tau^*\omega=-\omega$), and a
$U$-equivariant map $\Phi\colon M\to\u^*$ which is a moment map for
the action of $U$.  In addition we impose the following two
conditions:
\begin{align}
\Phi\bigl(\tau(m)\bigr)
&=-\sigma\bigl(\Phi(m)\bigr),\label{equation;anti}\\
\tau(um) &=\sigma(u)\tau(m)\label{equation;semi}
\end{align}
for all $m\in M$ and $u\in U$.
\end{definition}

Consider the semidirect product $U\rtimes\{\pm1\}$, where the action
of $\{\pm1\}$ on $U$ is given by the involution $\sigma$.  (By
$\{\pm1\}$ we mean a copy of $\Z/2$ written multiplicatively.)  Then
\eqref{equation;semi} says that $U\rtimes\{\pm1\}$ acts on $M$ and
that the action is \emph{sesquisymplectic} in the sense that
$u^*\omega=\eps(u)\omega$, where $\eps\colon
U\rtimes\{\pm1\}\to\{\pm1\}$ denotes the canonical homomorphism.

Properties \eqref{equation;anti} and \eqref{equation;semi} are not
independent.  For $\xi\in\u$ we denote by $\xi_M$ the induced vector
field on $M$ and by $\Phi^\xi$ the function on $M$ defined by
$\Phi^\xi(m)=\langle\Phi(m),\xi\rangle$.

\begin{lemma}\label{lemma;dependent}
If $U$ is connected\upn, then \eqref{equation;anti} implies
\eqref{equation;semi}.  Conversely\upn, if \eqref{equation;semi} is
satisfied\upn, we can shift $\Phi$ by a suitable constant to ensure
that \eqref{equation;anti} holds.
\end{lemma}

\begin{proof}
Using the moment map condition $d\Phi^\xi=\iota(\xi_M)\omega$ we
obtain from \eqref{equation;anti} that
\begin{equation}\label{equation;infisemi}
\tau_*(\xi_M)=\bigl(\sigma(\xi)\bigr)_M
\end{equation}
for all $\xi\in\u$.  This implies \eqref{equation;semi} if $U$ is
connected.  Conversely, if \eqref{equation;semi} holds then
\eqref{equation;infisemi} holds, whether or not $U$ is connected.
Therefore
\begin{multline*}
\langle d\Phi,\xi\rangle =\iota(\xi_M)\omega
=\iota\bigl(\tau_*\sigma(\xi)_M\bigr)\omega
=-\iota\bigl(\tau_*\sigma(\xi)_M\bigr)\tau^*\omega
=-\tau^*\bigl(\iota(\sigma(\xi)_M)\omega\bigr)\\
=-\tau^*\langle d\Phi,\sigma(\xi)\rangle =-\langle
d(\Phi\circ\tau),\sigma(\xi)\rangle =-\langle
d(\sigma\circ\Phi\circ\tau),\xi\rangle
\end{multline*}
for all $\xi$ and so $d(\sigma\circ\Phi\circ\tau)=-d\Phi$.
Furthermore, it is easy to check that $\sigma\circ\Phi\circ\tau$ is
$U$-equivariant.  It follows that
$$
\sigma\circ\Phi\circ\tau=-\Phi+c
$$
for some $U$-invariant $c\in\u^*$.  Multiplying by $\sigma$ on the
left and by $\tau$ on the right we get $\Phi
=-\sigma\circ\Phi\circ\tau+\sigma(c)$, and therefore $\sigma(c)=c$.
Putting $\Phi'=\Phi+c'$, with $c'=-c/2$, we see that $\Phi'$ is a
$U$-equivariant moment map satisfying
$$
\sigma\circ\Phi'\circ\tau =\sigma\circ(\Phi+c')\circ\tau
=\sigma\circ\Phi\circ\tau+c' =-\Phi+c+c' =-\Phi',
$$
i.e.\ $\Phi'\circ\tau=-\sigma\circ\Phi'$.
\end{proof}

Henceforth we assume $U$ to be connected, apart from a few exceptional
situations where for technical reasons we need to allow disconnected
groups.  These exceptions will be pointed out explicitly.

The focus of our attention is the fixed-point manifold $M^\tau$.  Let
$K$ be the identity component of the fixed-point group $U^\sigma$ and
let $\q\subset\u$ be the $-1$-eigenspace of $\sigma$.  Then
$\u=\k\oplus\q$, where $\k=\Lie(K)$.  Identifying $\k^*$ with the
annihilator of $\q$ and $\q^*$ with the annihilator of $\k$, we have
$\u^*=\k^*\oplus\q^*$.  Choose an almost complex structure $J$ on $M$
which is $U$-invariant and $\tau$-anti-invariant.  The existence of
such a $J$ is guaranteed by Lemma
\ref{lemma;almostcomplex}\eqref{item;complex} in Appendix
\ref{section;antisymplectic}.

\begin{proposition}\label{proposition;lagrangian}
\begin{enumerate}
\item\label{item;stable}
If $M^\tau$ is nonempty\upn, it is a $U^\sigma$-stable totally real
Lagrangian submanifold of $M$\upn, and $\Phi(M^\tau)$ is contained in
$\q^*$.
\item\label{item;connected}
Every connected component of $M^\tau$ is $K$-stable.
\item\label{item;orbit}
The $K$-orbit through every $m$ in $M^\tau$ is a connected component
of $Um\cap M^\tau$.
\end{enumerate}
\end{proposition}

\begin{proof}
The inclusion $\Phi(M^\tau)\subset\q^*$ is obvious from
\eqref{equation;anti} and the $U^\sigma$-stability of $M^\tau$ is
obvious from \eqref{equation;semi}.  The fact that $M^\tau$ is
Lagrangian and totally real follows from Lemma
\ref{lemma;almostcomplex}\eqref{item;lagrangian} (applied to the group
$\ca G=U\rtimes\{\pm1\}$).  This proves \eqref{item;stable}.

\eqref{item;connected} follows immediately from \eqref{item;stable}
and the fact that $K$ is connected.

Let $m\in M^\tau$.  It follows from \eqref{item;connected} that
$Km\subset Um\cap M^\tau$.  To prove \eqref{item;orbit} it suffices to
show that $Km$ is open in $Um\cap M^\tau$, or equivalently that the
tangent space to $Km$ at $m$ is equal to $T_m(Um)\cap T_mM^\tau$.  A
general tangent vector to $Um$ at $m$ is of the form
$\xi_M=\xi^+_M+\xi^-_M$ with $\xi^+\in\k$ and $\xi^-\in\q$.  Then
$\tau_*(\xi_M)=\sigma(\xi)_M=\xi^+_M-\xi^-_M$, so $\xi_M$ is tangent
to $M^\tau$ if and only if $\tau_*(\xi_M)=\xi_M$ if and only if
$\xi^-_M=0$, i.e.\ $\xi_M=\xi^+_M\in T_m(Km)$.
\end{proof}

\begin{example}[tori]\label{example;torus}
Let $U$ be a torus and $\sigma$ the involution which sends $u$ to
$u\inv$.  Then $U^\sigma$ is the subgroup of $2$-torsion elements of
$U$ and condition \eqref{equation;anti} boils down to $\Phi$ being
\emph{invariant} under $\tau$.  This is the setting of Duistermaat's
paper \cite{duistermaat;convexity}.
\end{example}

\begin{example}[connections]\label{example;riemann}
This infinite-dimensional example was suggested to us by C. Woodward.
Consider an oriented Riemann surface $\Sigma$ and an
orientation-reversing involution $f\colon\Sigma\to\Sigma$ without
fixed points.  Let $P$ be a trivial principal $U$-bundle over $\Sigma$
and $\eu A$ the space of flat connections on $P$.  According to Atiyah
and Bott $\eu A$ is a Hamiltonian $\eu U$-manifold, where $\eu
U=C^\infty(\Sigma,U)$ is the group of gauge transformations.
Involutions on $\eu A$ and $\eu U$ are given by pullback:
$\tau(A)=f^*A$ and $\sigma(u)=f^*u$.  Then $\tau$ is antisymplectic
because $f$ reverses the orientation and $\sigma$ is a group
homomorphism.  One checks easily that $\eu A$ is a Hamiltonian $(\eu
U,\sigma)$-manifold.  Moreover, the Lagrangian $\eu A^\tau$ is
isomorphic to the space of connections on the quotient bundle $P/f$
over the non-orientable quotient surface $\Sigma/f$, and $\eu
U^\sigma$ is the gauge group $C^\infty(\Sigma/f,U)$.
\end{example}

\begin{example}[cotangent bundles]\label{example;cotangent}
Let $S$ be any manifold with an involution $\tau$ and let $M=T^*S$ be
the cotangent bundle of $S$, equipped with its standard symplectic
form.  The cotangent lift of $\tau$ is a symplectic involution of $M$.
We turn it into an antisymplectic involution, also called $\tau$, by
``reversing momenta'', i.e.\ by multiplying by $-1$ in the cotangent
directions.  The fixed-point set $M^\tau$ is equal to the
\emph{conormal bundle} of $S^\tau$.  Suppose that $U$ acts on $S$ in
such a way that $\tau(us)=\sigma(u)\tau(s)$ for $s\in S$, and lift the
action to $M$ in the canonical way.  It is then a straightforward
exercise to check that $M$ is a Hamiltonian $(U,\sigma)$-manifold.
The moment map is given by $\Phi^\xi(s,p)=p(\xi_{S,s})$ for
$\xi\in\u$, $s\in S$ and $p\in T^*_sS$.  For instance, we can take $S$
to be the group $U$ with involution $\sigma$.  If we identify $M=T^*U$
with $U\times\u^*$ by means of left translations, then $\tau$ is given
by $\tau(u,\xi) =\bigl(\sigma(u),-\sigma(\xi)\bigr)$, and
$M^\tau=U^\sigma\times\q^*$.  On $M$ there are three different
Hamiltonian $U$-actions, namely the actions induced by left
multiplication, right multiplication, and conjugation.  Thus $M$ is a
Hamiltonian $(U,\sigma)$-manifold in three different ways.
\end{example}

\begin{example}[vector spaces]\label{example;linear}
In this example $U$ need not be connected.  Let $(V,\omega)$ be a
finite-dimensional symplectic vector space and let $\tau\in\GL(V)$ be
an antisymplectic involution.  Assume that $U$ acts linearly and
symplectically on $V$ and that \eqref{equation;semi} holds for all
$u\in U$ and $m\in V$.  Let $\Phi_V$ be the standard quadratic moment
map for the $U$-action given by
\begin{equation}\label{equation;quadratic}
\Phi_V^\xi(v)=\frac1{2}\omega(\xi v,v).
\end{equation}
Then the quadruple $(V,\omega,\tau,\Phi_V)$ is a Hamiltonian
$(U,\sigma)$-manifold, i.e.\ \eqref{equation;anti} holds
automatically.  The fixed-point set $V^\tau$ is clearly nonempty and
is in fact a $U^\sigma$-stable Lagrangian subspace.  Furthermore, by
Lemma \ref{lemma;complex} there exists a $U$-invariant and
$\tau$-anti-invariant compatible complex structure $J$ on $V$.  Then
we can define a real inner product $\beta$ on $V$ by
$\beta(v,w)=\omega(v,Jw)$ and a Hermitian inner product $h$ by
$h=\beta+i\omega$.  Clearly $h$ is $U$-invariant and
$h\bigl(\tau(v),\tau(w)\bigr)=h(w,v)$.  As a unitary $U^\sigma$-module
(but not as a $U$-module) $V$ is canonically isomorphic to the
complexification of $V^\tau$, and $\tau$ is complex conjugation
relative to the real form $V^\tau$.  Let $\ca A\colon U\to\U(V)$ be
the homomorphism defining the action of $U$ on $V$.  We call the data
$(V,\omega,\tau,J,\ca A)$ (often abbreviated to $V$) a \emph{unitary
$(U,\sigma)$-module}.  See Section \ref{section;quantization} for a
geometric construction and a classification of such modules.
\end{example}

\begin{remark}
Let $V$ be a unitary $(U,\sigma)$-module.  Consider the unitary group
$\U(V)$ and the standard involution on $\U(V)$ defined by
$\sigma_V(x)=(x\inv)^t$, where ``$t$'' denotes transpose with respect
to the Hermitian form $h$.  It follows from \eqref{equation;semi} that
the representation $\ca A$ intertwines the involutions on $U$ and
$\U(V)$: $\ca A\sigma=\sigma_V\ca A$.  In fact, this property can be
taken as an alternative definition of a unitary $(U,\sigma)$-module.
Here is yet another, purely complex algebraic characterization.
Define
\begin{equation}\label{equation;bilinear}
b(v,w)=h\bigl(v,\tau(w)\bigr).
\end{equation}
Then $b$ is a nondegenerate symmetric complex-bilinear form and is
``twisted'' invariant under $U$ in the sense that
$b\bigl(uv,\sigma(u)w\bigr)=b(v,w)$ for all $u\in U$ and $v$, $w\in
V$.  Conversely, a $U$-module furnished with such a form can be
endowed with the structure of a unitary $(U,\sigma)$-module such that
\eqref{equation;bilinear} holds.
\end{remark}

\begin{example}[coadjoint orbits]\label{example;orbit}
Let $\lambda\in\u^*$ and let $M=U\lambda$ be the coadjoint orbit
through $\lambda$, equipped with its standard symplectic form.  Assume
that $-\sigma(\lambda)\in M$, so that $-\sigma(M)=M$.  Then we can
define an involution $\tau$ on $M$ by putting $\tau(x)=-\sigma(x)$;
this is antisymplectic because $\sigma$ is a Lie algebra homomorphism.
The moment map is simply the inclusion of $M$ into $\u^*$, and
therefore \eqref{equation;anti} holds by construction.  Note that the
fixed-point set $M^\tau$ is equal to $M\cap\q^*$, so if it is nonempty
we may assume that $\lambda$ is in $\q^*$.  We call the coadjoint
orbit $M$ \emph{symmetric} if it is of the form $U\lambda$ with
$\lambda\in\q^*$.  We assert that if $\lambda\in\q^*$ then $M^\tau$ is
a single $K$-orbit,
\begin{equation}\label{equation;connect}
M^\tau=U^\sigma\lambda=K\lambda.
\end{equation}
To prove this we invoke Proposition \ref{proposition;lagrangian},
according to which it is enough to show that $M^\tau$ is connected.
Observe that the involution on $\u$ maps $[\u,\u]$ into itself.  Let
$\bar U$ be any connected Lie group with Lie algebra $[\u,\u]$ such
that the involution on $\u$ descends to an involution $\bar\sigma$ on
$\bar U$.  Then $\bar U$ acts on $\u$ and $\u^*$, and the orbits of
$\bar U$ are the same as those of $U$.  It follows from the proof of
Proposition 8.8(ii) in Helgason \cite[Ch.\
VII]{helgason;differential-geometry-lie-groups} that $M^\tau=\bar
U^{\bar\sigma}\lambda$ for any such $\bar U$.  We can take $\bar U$ to
be the universal cover of $[U,U]$, in which case the group $\bar
U^{\bar\sigma}$ is connected by Theorem 8.2 of \cite[Ch.\
VII]{helgason;differential-geometry-lie-groups}.  The upshot is that
$M^\tau=\bar U^{\bar\sigma}\lambda$ is connected, which concludes the
proof of \eqref{equation;connect}.  This example is continued in
\ref{example;kostant}, \ref{theorem;orbitquant},
\ref{example;realflag}, \ref{example;cross}, and in Section
\ref{section;application}.
\end{example}

\begin{example}[symplectic quotients]\label{example;quotient}
Let $M$ be a Hamiltonian $(U,\sigma)$-manifold and let $H$ be a Lie
group which acts on $M$ in a Hamiltonian fashion with momentum map
$\Psi\colon M\to\h^*$.  Assume that the action of $H$ commutes with
that of $U$.  Assume also that $\tau$ maps $\Psi\inv(0)$ into itself
and maps $H$-orbits to $H$-orbits, i.e.\ for all $h\in H$ and $m\in M$
there exists $h'\in H$ such that $\tau(hm)=h'\tau(m)$.  (A sufficient
condition for this to hold is that $H$ should be connected and that
for every $\eta\in\h$ there should exist $\eta'$ such that
$\tau^*(\Psi^\eta)=\Psi^{\eta'}$.)  Assume furthermore that $H$ acts
properly and freely on the fibre $\Psi\inv(0)$, so that the quotient
$M\qu H=\Psi\inv(0)/H$ is a well-defined symplectic manifold.  The
action of $U$ and the map $\tau$ descend to $M\qu H$ and it is easy to
see that $M\qu H$ is a Hamiltonian $(U,\sigma)$-manifold.
\end{example}

\begin{example}[projective varieties]\label{example;projective}
Let $V$ be a unitary $(U,\sigma)$-module as in Example
\ref{example;linear}.  Taking $H=S^1$ in Example
\ref{example;quotient} and letting $H$ act by complex scalar
multiplication on $V$, we conclude that the projective space $\P
V=V\qu S^1$ is a Hamiltonian $(U,\sigma)$-manifold in a natural way.
Likewise, every nonsingular complex subvariety of $\P V$ which is
stable under $U$ and $\tau$ inherits the structure of a Hamiltonian
$(U,\sigma)$-manifold.
\end{example}

\begin{example}[toric varieties]\label{example;toric}
Let $(U,\sigma)$ be as in Example \ref{example;torus}, let $\Delta$ be
a simple lattice polytope in $\u^*$, and let $M$ be the symplectic
toric variety (``Delzant space'') associated with $\Delta$.  As
explained e.g.\ in \cite{meinrenken-sjamaar;singular}, this space can
be obtained from $T^*U$ as a symplectic quotient with respect to an
appropriate subtorus of $U$ and therefore, by Examples
\ref{example;cotangent} and \ref{example;quotient}, has the structure
of a Hamiltonian $(U,\sigma)$-manifold.  In fact, $M$ is a complex
projective variety in a natural way and as such is defined over $\Z$
and hence over $\R$.  It is not hard to see that the antisymplectic
involution $\tau$ on $M$ is the same as complex conjugation.  The real
part $M^\tau$ of $M$ is nonempty and connected.  See e.g.\
\cite{jurkiewicz;torus-embeddings-polyhedra}, where the cohomology of
$M^\tau$ with $\Z/2$-coefficients is computed.
\end{example}

\begin{example}\label{example;restrictionproduct}
The product of a Hamiltonian $(U_1,\sigma_1)$-manifold and a
Hamiltonian $(U_2,\sigma_2)$-manifold is a Hamiltonian $(U_1\times
U_2,\sigma_1\times\sigma_2)$-manifold.  Let $H$ be a $\sigma$-stable
closed subgroup of $U$.  For simplicity we denote the restriction of
$\sigma$ to $H$ also by $\sigma$.  Then every Hamiltonian
$(U,\sigma)$-manifold inherits the structure of an
$(H,\sigma)$-manifold.  In particular, the product of two Hamiltonian
$(U,\sigma)$-manifolds is a Hamiltonian $(U,\sigma)$-manifold with
respect to the diagonal action.  Observe that these facts are also
true if the groups in question are not connected.
\end{example}

\begin{example}[duplication]\label{example;duplication}
Let $(M,\omega,\Phi)$ be any Hamiltonian $U$-manifold.  Let $M^-$ be
the same manifold as $M$, but furnished with the opposite symplectic
form $-\omega$.  Define $M\dup=M\sqcup M^-$, the disjoint union of $M$
and $M^-$, and define $\tau$ to be the antisymplectic involution which
interchanges the two components.  It has obviously no fixed points.
Extend the $U$-action and the moment map from $M$ to $M\dup$ by
putting $u\tau(m)=\tau\bigl(\sigma(u)m\bigr)$ and
$\Phi\bigl(\tau(m)\bigr)=-\sigma\bigl(\Phi(m)\bigr)$ for $m\in M$.
Then $M\dup$ is a Hamiltonian $(U,\sigma)$-manifold.  For instance, if
$M=U\lambda$ is a coadjoint orbit of $U$, then $M\dup\cong
U\lambda\sqcup-\sigma(U\lambda)$.
\end{example}

\begin{example}\label{example;linearduplication}
The process of duplication has linear analogue, which runs as follows.
Let $V$ be any unitary $U$-module with Hermitian metric
$h=\beta+i\omega$.  Let $V^-$ be the vector space conjugate to $V$
equipped with the Hermitian form $\bar h=\beta-i\omega$.  Put
$V\dup=V\oplus V^-$ and let $\tau$ be the map which interchanges the
two components.  It is obvious that $\tau$ is conjugate linear and
antisymplectic.  Turn $V^-$ into a $U$-module by putting
$u\tau(v)=\tau\bigl(\sigma(u)v\bigr)$ for $v\in V$; then
$(V\dup,\tau)$ becomes a unitary $(U,\sigma)$-module as in Example
\ref{example;linear}.  The fixed-point set of $\tau$ is the diagonal.
\end{example}

Here is a general criterion for the fixed-point set $M^\tau$ in a
Hamiltonian $(U,\sigma)$-manifold $(M,\omega,\tau,\Phi)$ to be
nonempty.  The integrality condition is satisfied for instance if $M$
is prequantizable in the sense of Definition \ref{definition;quantum}.

\begin{proposition}
Suppose that $M$ is a compact connected $(U,\sigma)$-manifold of
dimension $4n$ and that there exists a class $c\in H^2(M,\Z)$ such
that $\tau^*c=-c$ and the image of $c$ in $H^2(M,\R)$ is equal to the
cohomology class of $\omega$.  Then $M^\tau$ is nonempty and
$b^n|_{M^\tau}$ is nonzero\upn, where $b\in H^2(M,\Z/2)$ denotes the
reduction modulo $2$ of $c$.
\end{proposition}

\begin{proof}
Observe that $\tau^*b^n=b^n$ and that $b^{2n}$ is the mod~$2$
fundamental class of $M$.  The result now follows from Theorem 7.4 of
Bredon \cite[Ch. VII]{bredon;introduction-transformation}.
\end{proof}

We do not know of a good criterion for the fixed-point set to be
connected.  In Examples \ref{example;linear}, \ref{example;orbit} and
\ref{example;toric} the total space $M$ is simply connected and the
Lagrangian $M^\tau$ is connected.  The same is true for $M=T^*U$ (as
in Example \ref{example;cotangent}), provided that $U$ is simply
connected.  But here is an example where $M$ is simply connected and
$M^\tau$ is disconnected.  Let $S$ be the two-sphere with the
involution $\tau$ given by rotation through $\pi$ about the vertical
axis, and let $U=\SO(3)$ act on $S$ in the usual way.  Then
$\tau(us)=\sigma(u)\tau(s)$, where $\sigma$ is the inner involution of
$U$ defined by the matrix
$$
\begin{pmatrix}
-1&0&0\\0&-1&0\\0&0&1
\end{pmatrix}.
$$
Let $M=T^*S$ as in Example \ref{example;cotangent}.  Then $M^\tau$ is
the union of two Lagrangian planes.  However, this example is
noncompact.  This leads to the following question.

\begin{problem}
Is $M^\tau$ connected whenever $M$ is compact and simply
connected\upn?
\end{problem}

\section{The main result}\label{section;convex}

Let $(M,\omega,\tau,\Phi)$ be a Hamiltonian $(U,\sigma)$-manifold.
Let $T$ be any maximal torus of $U$ and $\t^*_+$ a (closed)
fundamental Weyl chamber in $\t^*$.  For any subset $X$ of $M$ we
define $\Delta(X)$ to be the intersection $\Phi(X)\cap\t^*_+$.  A
well-known theorem of Kirwan \cite{kirwan;convexity-III} (as
generalized by Hilgert, Neeb and Plank
\cite{hilgert-neeb-plank;symplectic-convexity;compositio} and Sjamaar
\cite{sjamaar;convexity}) says that if $M$ is connected and $\Phi$ is
proper, then $\Delta(M)$ is a convex polyhedron in $\t^*$.  The main
result of this paper is a comparable assertion about the Lagrangian
submanifold $M^\tau$ of $M$.

For this result to hold, the torus $T$ and the chamber $\t^*_+$ need
to be chosen in a suitable ``position'' relative to the involution
$\sigma$, namely as follows.  As before, let $K$ be the identity
component of $U^\sigma$ and $\q\subset\u$ the $-1$-eigenspace of the
involution $\sigma$.  Choose a maximal abelian subspace $\a$ of $\q$
and a maximal abelian subalgebra $\t$ of $\u$ which contains $\a$.
Let $A=\exp\a$ and $T=\exp\t$ be the corresponding subtori of $U$.
Let $R$ be the root system of $(U,T)$, $R'\subset R$ the subsystem of
roots which vanish on $\a$, and $R^\a$ the restricted root system.
(See Appendix \ref{section;symmetric}.)  Fix sets of positive roots
$R'_+$ in $R'$ and $R^\a_+$ in $R^\a$.  Define the set of positive
roots $R_+$ in $R$ to be the union of $R'_+$ and those roots in $R$
which restrict to elements of $R^\a_+$.  Denote the corresponding
fundamental Weyl chambers in $\t^*$ and $\a^*$ by $\a^*_+$ and
$\t^*_+$, respectively.  Then by Lemma \ref{lemma;chambers},
$\t^*_+\cap\q^*=\t^*_+\cap\a^*$ is equal to $\a^*_+$.

Now consider the $K$-stable subset $\Phi(M^\tau)$ of $\u^*$ and the
intersection
$$
\Delta(M^\tau)=\Phi(M^\tau)\cap\t^*_+.
$$
Since $\Phi(M^\tau)\subset\q^*$, $\Delta(M^\tau)$ is equal to
$\Phi(M^\tau)\cap\a^*_+$.  Our main result is the following convexity
theorem for $M^\tau$.

\begin{theorem}\label{theorem;convex}
Assume that $M$ is connected\upn, $M^\tau$ is nonempty\upn, and the
moment map $\Phi$ is proper.
\begin{enumerate}
\item\label{item;convex}
$\Delta(M^\tau)$ is equal to $\Delta(M)\cap\a^*$ and is therefore a
closed convex polyhedral subset of $\a^*$.
\item\label{item;vertex}
If $\lambda\in\a^*_+$ is a vertex of $\Delta(M^\tau)$ and
$m\in\Phi\inv(\lambda)\cap M^\tau$\upn, then
$A\subset[U_\lambda,U_\lambda]\,U_m$.  In particular\upn, the vertices
of $\Delta(M^\tau)$ which are contained in the relative interior of
$\t^*_+$ are images under $\Phi$ of $A$-fixed points in $M^\tau$.
\end{enumerate}
\end{theorem}

We call $\Delta(M^\tau)$ the \emph{moment polyhedron} of $M^\tau$.
From \eqref{item;convex} it follows easily that
$\Phi(M^\tau)=\Phi(M)\cap\q^*$.  Also, \eqref{item;convex} can be
interpreted as saying that if a complete list of inequalities for the
polyhedron $\Delta(M)$ is known, then the inequalities for
$\Delta(M^\tau)$ can be obtained simply by adding to the list a number
of equalities, namely those describing the subspace $\a^*$.  However,
as we shall see in Section \ref{section;application}, this is an
inefficient method for finding the inequalities for $\Delta(M^\tau)$.

In Sections \ref{section;kahler} and \ref{section;affine} we prove the
theorem for projective and affine varieties.  The general case is
treated in Section \ref{section;proof}, where we reduce it to the
algebraic case by means of a local normal form theorem.  The algebraic
versions of Theorem \ref{theorem;convex} enable us to generalize it to
various situations where $M$ is singular (Theorem
\ref{theorem;singular} and its corollaries), or where $\Phi$ is not
proper (Theorems \ref{theorem;affineconvex} and
\ref{theorem;cotangent}).  We also prove a result, Theorem
\ref{theorem;localcone}, which describes, for any $m\in M^\tau$
mapping into $\Delta(M^\tau)$, the shape of the moment polyhedron near
$\Phi(m)$ in terms of the coordinate ring of an affine variety
determined by infinitesimal data at $m$.

\begin{example}\label{example;identity}
If $\sigma$ is the identity map on $U$, then $\a=\{0\}$.  The
Lagrangian $M^\tau$ is $U$-invariant and satisfies
$\Phi(M^\tau)=\{0\}$, so Theorem \ref{theorem;convex} is trivially
true.
\end{example}

\begin{example}
If $M$ is compact and $U$ is a torus, Theorem \ref{theorem;convex}
says that $\Delta(M^\tau)$ is the convex hull of $\Phi(M^A\cap
M^\tau)$.  If in addition $\sigma(u)=u\inv$, then
$\Phi(M^\tau)=\Phi(M)$, which is the main result of Duistermaat's
paper \cite{duistermaat;convexity}.
\end{example}

\begin{example}\label{example;kostant}
Let $\lambda\in\a^*_+=\t^*_+\cap\q^*$ and let $M$ be the symmetric
coadjoint orbit $U\lambda$ of Example \ref{example;orbit}.  Then
$M^\tau=K\lambda$.  Consider the action of the maximal torus $T$ on
$M$.  By Theorem \ref{theorem;convex}, $\Delta(M^\tau) =\hull(M^A\cap
M^\tau)$.  According to Lemma \ref{lemma;weyl}\eqref{item;suborbit},
$M^A\cap M^\tau$ is equal to $W^\a\lambda$, the restricted Weyl group
orbit through $\lambda$.  The conclusion is that $\Delta(M^\tau)=\hull
W^\a\lambda$.  This was proved by Kostant in
\cite{kostant;convexity-weyl}.  Here is a slightly different way of
recovering the same result: Kostant also proved that $\Delta(M)=\hull
W\lambda$.  Therefore, by Theorem \ref{theorem;convex} and Lemma
\ref{lemma;weyl}\eqref{item;hull},
$$
\Delta(M^\tau) =\Delta(M)\cap\a^*=(\hull W\lambda)\cap\a^* =\hull
W^\a\lambda.
$$
See Section \ref{section;application} for yet another proof.
\end{example}

We emphasize that $M^\tau$ is not stable under the action of the torus
$A$.  Duistermaat has given a more intrinsic description of the set
$M^A\cap M^\tau$, which does not refer to the action of $A$ on the
ambient manifold $M$.  The submanifold $M^\tau$ is stable under the
$2$-torsion subgroup $A[2]=A^\sigma$, and clearly $M^A\cap M^\tau$ is
included in $(M^\tau)^{A[2]}$.  This inclusion is often, but not
always, an equality.  (For an example where it is not, let $U=\SU(2)$,
$\sigma=$ complex conjugation, so that $U^\sigma=K=\SO(2)$ and $A=T$.
Let $M$ be a coadjoint $U$-orbit.  Then $M$ is a sphere, $M^\tau$ is a
great circle, $M^A\cap M^\tau=M^A$ consists of the two poles, but
$A[2]$ is the centre of $U$, so it acts trivially on $M$.)  To remedy
this situation, define the sequence of submanifolds $Y_k$ of $M^\tau$
by $Y_k =M^{A[2^k]}\cap M^\tau$ for $k\ge0$, where $A[2^k]$ denotes
the $2^k$-torsion subgroup of $A$.  One easily proves that $A[2^k]$
maps $Y_{k-1}$ into itself and thus that $Y_k=Y_{k-1}^{A[2^k]}$ for
$k\ge1$.  Because $\{Y_k\}$ is a nested sequence of submanifolds, it
is eventually constant.  Since $\bigcup_{k\ge l}A[2^k]$ is dense in
$A$, we obtain $M^A\cap M^\tau=Y_l$ for sufficiently large $l$.  (Cf.\
\cite[Lemma 3.3]{duistermaat;convexity}.)

\section{Geometric quantization}\label{section;quantization}

The proof of Theorem \ref{theorem;convex} in the algebraic case makes
use of geometric quantization.  Given a Hamiltonian
$(U,\sigma)$-manifold $M$ which admits a K\"ahler polarization, how is
the involution on $M$ reflected in the quantization $Q(M)$?  In this
section we show that, subject to a compatibility condition on the
involution and the polarization, $Q(M)$ is a unitary
$(U,\sigma)$-module as in Example \ref{example;linear}.  The
compatibility condition is for example satisfied for a class of
coadjoint orbits which contains, but in some cases is strictly larger
than, the class of symmetric coadjoint orbits.  For such an orbit the
quantization is a unitary $(U,\sigma)$-module which is irreducible as
a $U$-module.  As a byproduct we obtain a characterization in terms of
the orbit method of representations of real and quaternionic type, and
a classification of unitary $(U,\sigma)$-modules.

The details are as follows.  Let $(M,\omega,\tau,\Phi)$ be a
Hamiltonian $(U,\sigma)$-manifold.  Suppose that the symplectic form
is integral.  Let $L$ be a prequantum line bundle, that is a complex
line bundle with Hermitian fibre metric $h$ and Hermitian connection
$\nabla$ with curvature $-2\pi i\omega$.

\begin{definition}\label{definition;quantum}
The prequantum bundle $L$ is \emph{$(U,\sigma)$-equivariant} if there
exist a lifting of the $U$-action on $M$ to a unitary action on $L$
and a lifting of $\tau$ to a real linear bundle map $\tau_L\colon L\to
L$ which is conjugate linear ($\tau_L(cl)=\bar c\tau_L(l)$ for all
$l\in L$ and $c\in\C$), maps horizontal subspaces to horizontal
subspaces, and is an involution ($\tau_L^2=\id_L$).
\end{definition}

The lifting $\tau_L$ defines an isomorphism $(L,\nabla)\to(\bar
L,\bar\nabla)$ of complex line bundles with connection and therefore
induces an an equivalence $(\tau^*L,\tau^*\nabla)\cong(\bar
L,\bar\nabla)$.

The vector fields generating the $U$-actions on $M$ and $L$ are
related by Kostant's formula
\begin{equation}\label{equation;lift}
\xi_L=(\xi_M)\hor+2\pi\Phi^\xi\nu_L
\end{equation}
for $\xi\in\u$.  (See e.g.\
\cite{guillemin-sternberg;geometric-quantization,%
kostant;quantization-prequantization}.)  Here ``hor'' refers to the
horizontal lift of a vector field and $\nu_L$ denotes the
infinitesimal generator of the scalar $S^1$-action on $L$.  This shows
that the obstruction to lifting the $U$-action is purely topological
and that the lift is unique, if it exists.

Let us now analyse the condition on $\tau$.  Because $\tau$ is
antisymplectic, the curvature of $(\tau^*L,\tau^*\nabla)$ is equal to
that of $(\bar L,\bar\nabla)$, and therefore $L'=\tau^*L\otimes\bar
L^*$ is a bundle with flat connection
$\nabla'=\tau^*\nabla\otimes1+1\otimes\bar\nabla^*$.  Hence, according
to the exact sequence
$$
H^1(M,\Z)\inj H^1(M,\R)\map H^1(M,S^1)\map H^2(M,\Z)\map
H^2(M,\R)\map\cdots,
$$
the pair $(L',\nabla')$ is represented up to equivalence by a class in
$H^1(M,S^1)$.  We conclude that a connection preserving conjugate
linear lifting $\tau_L$ of $\tau$ exists if and only if this class
vanishes.  According to Proposition 1.12.2 of Kostant
\cite{kostant;quantization-prequantization}, any another such lifting
$\tau_L'$ is given by $\tau_L'=c\tau_L$ for some nonzero complex
number $c$.  Moreover, we can choose $c$ such that $\tau_L'$ maps $h$
to $\bar h$, that is to say, preserves the real inner product $\Re h$.
Such a lifting is ``almost'' an involution.

\begin{lemma}\label{lemma;lift}
Let $\tau_L$ be conjugate linear lifting of $\tau$ which preserves the
connection.
\begin{enumerate}
\item\label{item;inner}
If $\tau_L$ preserves $\Re h$\upn, then $\tau_L^2=a\,\id_L$ for some
complex number $a$ of norm $1$.
\item\label{item;involution} 
If $\tau_L$ is an involution\upn, then $\tau_L$ preserves $\Re h$.
\item\label{item;phase}
Let $\tau_L'$ be another conjugate linear lifting of $\tau$ which
preserves the connection.  If $\tau_L$ is an involution\upn, then
$\tau_L'$ is an involution if and only if it is of the form
$\tau_L'=c\tau_L$ with $\lvert c\rvert=1$.
\end{enumerate}
\end{lemma}

\begin{proof}
Assume that $\tau_L$ preserves $\Re h$.  Since $\tau$ is an
involution, the map $\tau_L^2$ is then a (complex linear) automorphism
of $L$ which covers the identity and preserves $\nabla$ and $h$.  It
follows from the aforementioned result of Kostant's that
$\tau_L^2=a\,\id_L$ for some $a$ in the unit circle.  This proves
\eqref{item;inner}.

Again by Kostant's result, the fact that $\tau_L$ is conjugate linear
implies that $\tau_L^*h=b\bar h$ for some $b>0$, and hence
$\tau_L^*\Re h=b\Re h$.  If $\tau_L$ is an involution, then $\Re
h=(\tau_L^*)^2\Re h=b^2\Re h$, so $b=1$ and $\tau_L^*\Re h=\Re h$.
This proves \eqref{item;involution}.

Finally, assume that $\tau_L$ is an involution.  We know that
$\tau_L'=c\tau_L$ for some $c\ne0$.  Therefore $\tau_L'$ is an
involution if and only if $\id_L=(\tau_L')^2=c\tau_Lc\tau_L=c\bar
c\tau_L^2=\lvert c\rvert^2\id_L$, i.e.\ $\lvert c\rvert=1$.  This
proves \eqref{item;phase}.
\end{proof}

\begin{remark}
Even if $\tau$ is liftable, it is not always possible to find an
involutive lift.  In some cases one can only find a lift $\tau_L$
which satisfies the ``quaternionic'' condition $\tau_L^2(l)=-l$.  See
also Remark \ref{remark;quaternion}.
\end{remark}

Suppose from now on that $M$ admits an \emph{integrable} complex
structure $J$ which is $U$-invariant, $\tau$-anti-invariant and
compatible with $\omega$.  Then $\omega$ is a positive $(1,1)$-form,
so $L$ admits a holomorphic structure such that the projection $L\to
M$ is holomorphic.  In addition, we can choose the metric connection
$(\nabla,h)$ to be compatible with the holomorphic structure.

\begin{lemma}\label{lemma;connection}
Let $\tau_L\colon L\to L$ be any real linear bundle map which lifts
$\tau$ and preserves the real inner product $\Re h$.  Then $\tau_L$ is
antiholomorphic if and only if $\tau_L$ is conjugate linear and
preserves the connection.
\end{lemma}

\begin{proof}
Assume that $\tau_L$ is antiholomorphic.  Then the derivative of
$\tau_L$ at any point in $L$ is a conjugate linear map.  By taking a
point $m$ in the zero section we see that the restriction of
$(\tau_L)_*$ to $L_m$ is a conjugate linear map $L_m\to L_{\tau(m)}$,
which implies that the restriction of $\tau_L$ to every fibre is
conjugate linear.  Let $f\colon L\to\R$ be the smooth function defined
by $f(l)=h(l,l)$.  The fact that the connection is Hermitian implies
that the horizontal subspace at $l$ is equal to
$$
\ca H_l=\ker d_lf\cap i\ker d_lf.
$$
Together with the fact that $\tau_L$ is antiholomorphic and preserves
$f$ this implies that $(\tau_L)_*\ca H_l=\ca H_{\tau_L(l)}$ for all
$l\in L$, that is to say $\tau_L$ preserves the connection.

Now assume that $\tau_L$ is conjugate linear and preserves the
connection.  The tangent space to $L$ at any point $l$ decomposes as
$T_lL=\ca H_l\oplus L_m$, where $m$ is the base point of $l$.  Because
$\tau_L$ preserves the connection, its derivative at $l$ is the direct
sum of two maps $L_m\to L_{\tau(m)}$ and $\ca H_l\to\ca
H_{\tau_L(l)}$.  Both of these maps are conjugate linear, the first
one by assumption and the second one because $\tau$ is antiholomorphic
and the connection is complex.  Therefore $(\tau_L)_*$ is conjugate
linear, i.e.\ $\tau_L$ is antiholomorphic.
\end{proof}

Henceforth we assume $L$ to be $(U,\sigma)$-equivariant and fix an
involutive lifting $\tau_L$ as in Definition \ref{definition;quantum}.
Then $\tau_L$ preserves $\Re h$ by Lemma \ref{lemma;lift} and hence is
antiholomorphic by Lemma \ref{lemma;connection}.  Let
$Q(M)=\Gamma(M,L)$ be the space of global holomorphic sections of $L$.
We define an involution on the global \emph{smooth} sections of $L$ by
sending a section $s$ to $\tau_L\circ s\circ\tau$.  Because $\tau$ and
$\tau_L$ are both antiholomorphic, this involution maps $Q(M)$ into
itself.  We denote its restriction to $Q(M)$ by $\tau_Q$.  Since
$\tau_L$ is conjugate linear, so is $\tau_Q$.  Moreover, $\tau_Q$
preserves the real part of the Hermitian inner product on $Q(M)$ given
by $\langle s,t\rangle=\int_Mh\bigl(s(m),t(m)\bigr)\,dm$, where
$dm=\lvert\omega^n/n!\rvert$ is Liouville measure on $M$.  With this
inner product, $Q(M)$ is also a unitary $U$-module in a natural way.
The involution and the $U$-representation on $Q(M)$ are related in the
expected manner.  Since $\tau_L$ preserves the connection, we have
$(\tau_L)_*X\hor=(\tau_*X)\hor$ for all tangent vectors $X$.  Together
with \eqref{equation;anti}, \eqref{equation;infisemi} and
\eqref{equation;lift} this implies that for all $\xi\in\u$
\begin{equation}\label{equation;horizontal}
\begin{split}
(\tau_L)_*\xi_L &=(\tau_L)_*\bigl((\xi_M)\hor\bigr)
+2\pi(\Phi^\xi\circ\tau)(\tau_L)_*\nu_L\\
&=(\tau_*\xi_M)\hor-2\pi\Phi^{\sigma(\xi)}(-\nu_L)\\
&=\bigl(\sigma(\xi)_M\bigr)\hor+2\pi\Phi^{\sigma(\xi)}\nu_L\\
&=\sigma(\xi)_L.
\end{split}
\end{equation}
Therefore $\tau_L(ul)=\sigma(u)\tau_L(l)$ for all $u\in U$ and $l\in
L$.  From this we obtain easily that $\tau_Q(us)=\sigma(u)\tau_Q(s)$
for all $u\in U$ and holomorphic sections $s$.  The result is as
follows.

\begin{proposition}\label{proposition;quant}
The vector space $Q(M)$ is a unitary $(U,\sigma)$-module as in Example
{\rm\ref{example;linear}}\upn, finite-dimensional if $M$ is compact.
Two different involutive connection-preserving liftings of $\tau$ give
rise to isomorphic unitary $(U,\sigma)$-modules.
\end{proposition}

\begin{proof}
It remains to prove the uniqueness part.  According to Lemma
\ref{lemma;lift}, any other lifting $\tau_L'$ with the stated
properties is of the form $\tau_L'=c\tau_L$ with $\lvert c\rvert=1$.
Then $\tau_L'$ gives rise to the involution $\tau_Q'=c\tau_Q$ of
$Q(M)$.  The scalar multiplication operator $C(s)=\sqrt{c}\,s$ is
unitary, commutes with $U$ and satisfies $C\tau_QC\inv=C^2\tau_Q
=\tau_Q'$, so defines a unitary intertwining operator between the
$(U,\sigma)$-modules $\bigl(Q(M),\tau_Q\bigr)$ and
$\bigl(Q(M),\tau_Q'\bigr)$.
\end{proof}

Kodaira's Embedding Theorem implies that if $M$ is compact, it is
biholomorphically equivalent to a complex projective variety.  Using
Proposition \ref{proposition;quant} one can easily show that Kodaira's
map, $j$, is $U\rtimes\{\pm1\}$-equivariant, and hence that its image
$j(M)$ is a $U$-stable and $\tau$-stable nonsingular projective
variety, as in Example \ref{example;projective}.  In general $j$ is
not a symplectomorphism, but a straightforward Darboux-type argument
shows that there exists a $U\rtimes\{\pm1\}$-equivariant
diffeomorphism $f$ from $M$ to itself such that $j\circ f$ is a
symplectomorphism.  The conclusion is as follows.

\begin{corollary}\label{corollary;projective}
Every compact Hamiltonian $(U,\sigma)$-manifold admitting a
$(U,\sigma)$-equivariant K\"ahler prequantization is symplectically
and $(U,\sigma)$-equivariantly diffeomorphic to a complex projective
$(U,\sigma)$-manifold.
\qed
\end{corollary}

Using Proposition \ref{proposition;quant} we can give a complete
description of all finite-dimensional $(U,\sigma)$-modules.  Let us
first investigate the irreducible case.  It is well-known that the
K\"ahler quantization of a Hamiltonian $U$-manifold is an irreducible
$U$-module if and only if the manifold is isomorphic to an integral
coadjoint orbit, i.e.\ an orbit of the form $U\lambda$ with $\lambda$
integral dominant.  Such an orbit is diffeomorpic in a natural way to
the complex homogeneous space $U^\C/P_\lambda$, where $U^\C$ is the
complexification of $U$ and $P_\lambda$ is the parabolic subgroup
associated with $\lambda$.  It follows from the Borel-Weil Theorem
that the quantization of $U\lambda$ is the irreducible module
$V_\lambda$ with highest weight $\lambda$.  We shall establish a
necessary and sufficient criterion for $V_\lambda$ to be a
$(U,\sigma)$-module.

First we need an auxiliary result.  Let $K'$ be the centralizer of $A$
in $K$.  Then $T'=(K')_0\cap T$ is a maximal torus of $K'$.  Let
$W'\subset W$ the Weyl group of (the identity component of) $K'$ with
respect to $T'$ and let $w'_0$ be the longest element of $W'$.  Define
$\sigma_+\colon\t^*\to\t^*$ by $\sigma_+(\lambda)
=-\sigma(w'_0\lambda)$.  As before, $R_+$ denotes a set of positive
roots in $\t^*$ chosen as in \eqref{equation;positive}.

\begin{lemma}\label{lemma;centrallong}
\begin{enumerate}
\item\label{item;chamber}
$\sigma_+$ is an involution and $\sigma_+(R_+)=R_+$.  Hence
$\sigma_+(\t^*_+)=\t^*_+$ and $\sigma_+$ preserves the partial
ordering on the weight lattice defined by $R_+$.
\item\label{item;fix}
$\sigma_+(\lambda)=\lambda$ if and only if either
$\lambda\in\a^*$\upn, or $\lambda\in(\t')^*$ and
$\lambda=-w'_0\lambda$.
\item\label{item;representative}
Let $k_0\in\ca N_{(K')_0}(T')$ be any representative of $w'_0$.  Let
$\lambda$ be an integral weight and let $\chi_\lambda$ denote the
character of $T$ obtained by exponentiating $\lambda$.  If
$\lambda\in\a^*$\upn, then $\chi_\lambda(k_0^2)=1$.  If
$\lambda\in(\t')^*$ and $\lambda=-w'_0\lambda$\upn, then
$\chi_\lambda(k_0^2)=\pm1$.
\end{enumerate}
\end{lemma}

\begin{proof}
The action of $W'$ fixes $\a^*$ and therefore commutes with $\sigma$.
Since $w'_0$ and $\sigma$ are involutions, so is $\sigma_+$.  Observe
that $w_0'$ maps $R'_+$ to $R'_-$ and permutes $R_+\setminus R'_+$,
whereas $\sigma$ fixes $R'_+$ and sends $R_+\setminus R'_+$ to
$R_-\setminus R'_-$.  This implies that $w'_0\circ\sigma$ maps $R_+$
to $R_-$.  Therefore $\sigma_+(R_+)=R_+$.  This proves
\eqref{item;chamber}.

The fact that $w'_0$ fixes $\a^*$ implies that $\sigma_+$ respects the
decomposition \eqref{equation;cartandecompose}, and that
$\sigma_+|_{\a^*}=\id$ and $\sigma_+|_{(\t')^*}=-w'_0$.  This
immediately implies \eqref{item;fix}.

Let $\lambda$ be an integral weight.  Since $w'_0$ is an involution,
$k_0^2$ is in $T'$.  If $\lambda\in\a^*$, then $T'$ is contained in
the kernel of $\chi_\lambda$, so $\chi_\lambda(k_0^2)=1$.  Now suppose
$\lambda\in(\t')^*$ satisfies $\lambda=-w'_0\lambda$.  Then
$\chi_\lambda(k_0tk_0\inv) =\chi_\lambda(t)\inv$ for all $t\in T$.
Taking $t=k_0^2$ we obtain $\chi_\lambda(k_0^2)
=\chi_\lambda(k_0^2)\inv$.
\end{proof}

\begin{theorem}\label{theorem;orbitquant}
Let $\lambda$ be an integral dominant weight of $U$.  Then the
irreducible $U$-module $V_\lambda$ is a $(U,\sigma)$-module if and
only if $\sigma_+(\lambda)=\lambda$ and $\chi_\lambda(k_0^2)=1$\upn,
where $k_0$ is as in Lemma
{\rm\ref{lemma;centrallong}\eqref{item;representative}}.
\end{theorem}

In particular, $V_\lambda$ is a $(U,\sigma)$-module if
$\lambda\in\a^*$.  It is easy to check that if $\sigma_+(\lambda)
=\lambda$ the condition $\chi_\lambda(k_0^2)=1$ does not depend on the
choice of the representative $k_0$ of $w'_0$.

\begin{proof}
Let $(V,\tau)$ be any $(U,\sigma)$-module and let $v$ be a weight
vector in $V$ with weight $\mu$.  The rule $\tau(uv)=\sigma(u)\tau(v)$
implies that $\tau(v)$ is a weight vector with weight $-\sigma(\mu)$,
and therefore
\begin{equation}\label{equation;weightvector}
\text{$k_0\tau(v)$ is a weight vector with weight $\sigma_+(\mu)$.}
\end{equation}
Now consider the orbit $M=U\lambda$ and the $U$-module
$Q(M)=V_\lambda$.  Assume that $V_\lambda$ is a $(U,\sigma)$-module.
Let $v$ be a highest-weight vector.  Then $k_0\tau(v)$ has weight
$\sigma_+(\lambda)$ according to \eqref{equation;weightvector}.  We
conclude from Lemma \ref{lemma;centrallong}\eqref{item;chamber} that
$\sigma_+(\lambda)=\lambda$ and $k_0\tau(v)=cv$ for some $c\in\C$.
This implies $(k_0\tau)^2(v)=\lvert c\rvert^2v$ because $\tau$ is
conjugate linear.  Hence
$$
\chi_\lambda(k_0^2)\,v=k_0^2v=k_0^2\tau^2(v)=(k_0\tau)^2(v)=\lvert
c\rvert^2v,
$$
where we have used that $\tau$ commutes with $k_0$ since $k_0\in K$.
This implies $\chi_\lambda(k_0^2)=1$, since
$\lvert\chi_\lambda(k_0^2)\rvert=1$.

Now assume that $\sigma_+(\lambda)=\lambda$ and
$\chi_\lambda(k_0^2)=1$.  Then $-\sigma(\lambda)=w'_0\lambda$, so by
Example \ref{example;orbit} the orbit $M$ has a well-defined
antisymplectic involution $\tau=-\sigma$.  Let us check that $\tau$ is
antiholomorphic.  Because the $U$-action is holomorphic and
transitive, we need only check that $(\ca L_{k_0})_*\circ\tau_*\colon
T_\lambda M\to T_\lambda M$ is conjugate linear, where ``$\ca L$''
denotes left multiplication on $M$ by elements of $U$.  Since $w'_0$
is an involution,
\begin{equation}\label{equation;orbitinvolution}
\tau(u\lambda)=\sigma(u)\tau(\lambda)=\sigma(u)w'_0\lambda
=\sigma(u)k_0\inv\lambda.
\end{equation}
Let us identify $T_\lambda M$ with $\u/\u_\lambda$.  Left-multiplying
by $k_0$ and differentiating at $1$ we obtain from
\eqref{equation;orbitinvolution} the commutative diagram
\begin{equation}\label{equation;adsigma}
\begin{CD}
\u@>\Ad(k_0)\circ\sigma>>\u\\
@VVV@VVV\\
\u/\u_\lambda@>{(\ca L_{k_0})_*\circ\tau_*}>>\u/\u_\lambda.
\end{CD}
\end{equation}
Consider the root-space decomposition
$\u^\C=\t^\C\oplus\bigoplus_{\alpha\in R}\u^\C_\alpha$ and select root
vectors $E_\alpha\in\u^\C_\alpha\setminus\{0\}$ such that
$E_{-\alpha}=-\bar E_\alpha$ for all $\alpha\in R$.  Denote the
complex linear extension of $\Ad(k_0)\circ\sigma$ by $C$.  Observe
that $C\u^\C_\alpha=\u^\C_{-\sigma_+(\alpha)}$, which implies that
there exist $c_\alpha\in\C$ such that
\begin{equation}\label{equation;constants}
CE_\alpha=c_\alpha E_{-\sigma_+(\alpha)}\qquad\text{and}\qquad
CE_{-\alpha}=\bar c_\alpha E_{\sigma_+(\alpha)}
\end{equation}
for $\alpha\in R_+$.  Let $R^\lambda=R_+^\lambda\cup R_-^\lambda$ be
the set of all roots which annihilate $\lambda$.  That is to say,
$\alpha\in R^\lambda$ if and only $\alpha\spcheck(\lambda)=0$, where
$\alpha\spcheck$ denotes the dual root of $\alpha$.  Then
$\sigma_+(R_\pm^\lambda) =R_\pm^\lambda$ because
$\sigma_+(\lambda)=\lambda$ and $\sigma_+$ maps $R_+$ into itself.
Hence the Lie algebra of the parabolic subgroup $P_\lambda$ is given
by
$$
\p_\lambda=\t^\C\oplus\bigoplus_{\alpha\in R_+}\u^\C_\alpha
\oplus\bigoplus_{\alpha\in R_-^\lambda}\u^\C_\alpha.
$$
Therefore
\begin{equation}\label{equation;tangent}
\begin{split}
\u^\C/\p_\lambda &\cong\bigoplus_{\alpha\in R_-\setminus
R_-^\lambda}\C E_\alpha,\\
T_\lambda M\cong\u/\u_\lambda &\cong\bigoplus_{\alpha\in R_-\setminus
R_-^\lambda}(\R X_\alpha\oplus\R Y_\alpha),
\end{split}
\end{equation}
where $X_\alpha=\frac1{2i}(E_{-\alpha}+E_\alpha)$ and
$Y_\alpha=\frac1{2}(E_{-\alpha}-E_\alpha)$.  The complex structure on
$T_\lambda M$ is by definition the pullback of the complex structure
on $\u^\C/\p_\lambda$, so we obtain from \eqref{equation;tangent} that
$JX_\alpha=Y_\alpha$ and $JY_\alpha=-X_\alpha$.  It follows from
\eqref{equation;constants} that for all $\alpha\in R_-\setminus
R_-^\lambda$
\begin{align*}
CX_\alpha &=a_\alpha X_{\sigma_+(\alpha)}+b_\alpha
Y_{\sigma_+(\alpha)},\\
CY_\alpha &=b_\alpha X_{\sigma_+(\alpha)}-a_\alpha
Y_{\sigma_+(\alpha)},
\end{align*}
where $a_\alpha=\Re c_\alpha$ and $b_\alpha=\Im c_\alpha$.  Thus
$CJ=-JC$.  Together with \eqref{equation;adsigma} this proves that
$(\ca L_{k_0})_*\circ\tau_*$ is conjugate linear on $\u/\u_\lambda$.

Next we show that $U\lambda$ possesses a $(U,\sigma)$-equivariant
prequantum line bundle.  Let $L$ be the homogeneous bundle
$$
L=U\times^{U_\lambda}\C_\lambda\cong U^\C\times^{P_\lambda}\C_\lambda,
$$
where $\C_\lambda$ denotes the obvious one-dimensional representation
of $U_\lambda$, resp.\ $P_\lambda$.  This bundle has a $U$-invariant
Hermitian structure, namely the one induced by the standard inner
product on $\C$.  Put
$$
\tau_L[u,z]=[\sigma(u)k_0\inv,\bar z].
$$
The assumption that $\sigma_+(\lambda)=\lambda$ implies that $\tau_L$
is well-defined as a map from $L$ to itself.  The assumption that
$\chi_\lambda(k_0^2)=1$ implies that it is an involution.  It is clear
that $\tau_L$ is a conjugate linear bundle map and it follows from
\eqref{equation;orbitinvolution} that it is a lifting of the
involution $\tau$ on $M$.  The length function on $L$ is given by
$f([u,z])=\lvert z\rvert^2$, which is evidently preserved by $\tau_L$.
One checks directly that $\tau_L(ul)=\sigma(u)\tau_L(l)$, so
$(\tau_L)_*\xi_L=\sigma(\xi)_L$ for all $\xi\in\u$.  Furthermore,
$(\tau_L)_*\nu_L=-\nu_L$ and $\tau^*\Phi^\xi=-\Phi^{\sigma(\xi)}$, so
by reading \eqref{equation;horizontal} backwards we get
$(\tau_L)_*\bigl((\xi_M)\hor\bigr)=(\tau_*\xi_M)\hor$ for all $\xi$.
Since $M$ is homogeneous, this implies $(\tau_L)_*\ca H_l=\ca
H_{\tau_L(l)}$ for all $l\in L$, and so $\tau_L$ preserves the
connection.  Therefore $\tau_L$ is antiholomorphic by Lemma
\ref{lemma;connection}.  The upshot is that $L$ is a
$(U,\sigma)$-equivariant prequantum line bundle on $M$.  Hence, by
Proposition \ref{proposition;quant}, $V_\lambda=Q(M)$ is a unitary
$(U,\sigma)$-module.
\end{proof}

\begin{remark}[quaternionic modules]\label{remark;quaternion}
If $\sigma_+(\lambda)=\lambda$ but $\chi_\lambda(k_0^2)=-1$, then the
lifting $\tau_L$ is not an involution, but satisfies
$\tau_L^2[u,z]=[u,-z]$.  Otherwise the proof goes through as stated
and the result is a unitary irreducible $U$-module $Q(M)$ with an
antiholomorphic orthogonal map $\tau_Q$ satisfying $\tau_Q^2=-\id$ and
$\tau_Q(uv)=\sigma(u)v$.  In the extreme case where $\sigma$ is the
identity such a module is known as a $U$-module of \emph{quaternionic}
(or \emph{symplectic}) type, and what we have called a
$(U,\sigma)$-module is nothing but a $U$-module of \emph{real} (or
\emph{orthogonal}) type.  The condition $\sigma_+(\lambda)=\lambda$ is
then equivalent to $\lambda=-w_0\lambda$ ($w_0$ being the longest
element in $W$), in other words to the representation $V_\lambda$
being self-conjugate.
\end{remark}

\begin{example}\label{example;sun}
If $U=\SU(n)$ and $\sigma$ is complex conjugation, then
$U^\sigma=K=\SO(n)$ and $\a=\t$.  Hence \emph{every} coadjoint orbit
is symmetric and every irreducible $U$-module admits the structure of
a $(U,\sigma)$-module.  On the fundamental representation
$\bigwedge^k\C^n$ the involution is simply complex conjugation with
respect to a real form $\bigwedge^k\R^n$.  On the other hand, if
$\sigma$ is the identity, then $\a=\{0\}$ and only a few orbits give
rise to representations of real type.  For instance, $\bigwedge^k\C^n$
is self-dual if and only if $n$ is even and $k=n/2$.  It is of real
type if and only if $n\equiv0\bmod4$ and $k=n/2$.  (See e.g.\ Tits
\cite{tits;tabellen}.)
\end{example}

We can now easily describe the structure of an arbitrary
$(U,\sigma)$-module.

\begin{addendum}\label{addendum;classification}
\begin{enumerate}
\item\label{item;direct}
Every finite-dimensional object in the category of
$(U,\sigma)$-modules is isomorphic to a direct sum of irreducible
objects.
\item\label{item;irreducible}
The irreducible objects are
\begin{alignat}{2}
\label{equation;irreducible}
&V_\lambda&\quad&\text{with $\sigma_+(\lambda)=\lambda$ and
$\chi_\lambda(k_0^2)=1$},\\
\label{equation;dup}
&V_\lambda\dup&&\text{with $\sigma_+(\lambda)\ne\lambda$ or
$\chi_\lambda(k_0^2)=-1$}.
\end{alignat}
Here $\lambda$ ranges over the dominant weights of $U$ and
$V_\lambda\dup$ stands for the duplicated module defined in Example
{\rm\ref{example;linearduplication}}.
\item\label{item;dup}
As a $U$-module $V_\lambda\dup$ is isomorphic to $V_\lambda\oplus
V_{\sigma_+(\lambda)}$.
\end{enumerate}
\end{addendum}

Note that a module of type \eqref{equation;dup} can be described as
the quantization of the disconnected manifold $(U\lambda)\dup\cong
U\lambda\sqcup-\sigma(U\lambda)$.

\begin{proof}
\eqref{item;direct} is proved by induction on the dimension and
\eqref{item;dup} follows from \eqref{equation;weightvector} and Lemma
\ref{lemma;centrallong}\eqref{item;chamber}.

It is clear that the objects \eqref{equation;irreducible} are
irreducible and it follows from Theorem \ref{theorem;orbitquant} that
the objects \eqref{equation;dup} are irreducible as well.  Now let
$(V,\tau)$ be an arbitrary irreducible object.  If it is irreducible
as a $U$-module, we conclude from Theorem \ref{theorem;orbitquant}
that it must be of type \eqref{equation;irreducible}.  If not, it must
be of the form $V_\lambda\oplus\tau(V_\lambda)\cong V_\lambda\oplus
V_{\sigma_+(\lambda)}\cong V_\lambda\dup$ for some
$\lambda\in\Lambda^*_+$.  In this case we cannot have
$\sigma_+(\lambda)=\lambda$ and $\chi_\lambda(k_0^2)=1$, for then
$V_\lambda\dup$ would be isomorphic to $V_\lambda\oplus V_\lambda$ as
a $(U,\sigma)$-module.  (The easiest way to see this is by using the
bilinear form $b$ defined in \eqref{equation;bilinear}.  It is not
hard to show that if a $U$-module $E$ carries two nondegenerate
symmetric bilinear forms $b_1$, $b_2$ satisfying
$b_i\bigl(uv,\sigma(u)w\bigr)=b_i(v,w)$, then $f^*b_1=b_2$ for some
$U$-equivariant linear map $f$.  Applying this observation to
$E=V_\lambda\dup$ with $b_1$ being the form defined by the involution
$\tau$ and $b_2$ the one defined by the involutions on each of the
summands $V_\lambda$, we obtain $V_\lambda\dup\cong V_\lambda\oplus
V_\lambda$ as a $(U,\sigma)$-module.)  Thus we are in case
\eqref{equation;dup}.  This proves \eqref{item;irreducible}.
\end{proof}

\begin{remark}\label{remark;disconnect}
If $U$ is not connected, \eqref{item;direct} is still true, but the
description of the irreducible objects is less clear-cut.  However,
the fact remains that each such object is either of the form $V$ or
$V\dup$, where $V$ is irreducible as a $U$-module.
\end{remark}

\section{The projective case}\label{section;kahler}

In this section we prove Theorem \ref{theorem;convex} for compact
integral K\"ahler manifolds.  The gist of the argument is that the
Lagrangian $M^\tau$ is dense in $M$ for the complex Zariski topology.
Noteworthy is that the proof applies to singular as well as
nonsingular varieties, and that it leads to a stronger conclusion,
Corollaries \ref{corollary;orbitclosure} and \ref{corollary;dual},
which involve the dual symmetric pair $(G,K)$.  On top of this we
obtain a simple characterization of the polytope $\Delta(M^\tau)$ in
terms of certain rational weights associated with the quantization of
$M$, analogous to the characterization of $\Delta(M)$ given by
Guillemin and Sternberg \cite{guillemin-sternberg;convexity;;1982} and
Mumford \cite{ness;stratification}.

Let $(M,\omega,\tau,\Phi)$ be a compact connected Hamiltonian
$(U,\sigma)$-manifold equipped with a complex structure which is
compatible with $\omega$ and relative to which $\tau$ is
antiholomorphic.  Furthermore, as in Proposition
\ref{proposition;quant}, let $L$ be a holomorphic line bundle with a
Hermitian connection $\nabla$ such that $\curv\nabla=-2\pi i\omega$
and $\tau$ lifts to an involutive antiholomorphic bundle map $\tau_L$
which preserves $\nabla$.

According to Corollary \ref{corollary;projective} we can alternatively
think of $M$ as a $U$-stable and $\tau$-stable complex submanifold of
$\P V$ for some unitary $(U,\sigma)$-module $V$.

Following Brion \cite{brion;image}, we start by defining a
``rational'' analogue of the sets $\Delta(M)$ and $\Delta(M^\tau)$.
Let $\Lambda^*_+$ be the monoid of dominant weights in $\t^*_+$.

\begin{definition}\label{definition;brion}
The \emph{highest-weight set} of $M$ is the subset $\ca C(M)$ of
$\Lambda^*\otimes\Q$ consisting of all quotients $\lambda/n$, where
$\lambda\in\Lambda^*_+$ and $n$ is a positive integer such that the
irreducible representation $V_\lambda$ occurs in the $U$-module
$\Gamma(M,L^n)$.
\end{definition}

Since
$$
\Hom\bigl(\Gamma(M,L^n),V_\lambda\bigr)\cong
\Gamma(M,L^n)\otimes\Gamma(U\lambda^*,L_{\lambda^*})
\cong\Gamma(M\times U\lambda^*,L^n\boxtimes L_{\lambda^*}),
$$
we see that $\xi\in\ca C(M)$ if and only if $\xi=\lambda/n$ for some
$\lambda\in\Lambda^*_+$ and $n>0$ such that there exist $m\in M$,
$u\in U$, and a $U$-invariant global holomorphic section $s$ of
$L^n\boxtimes L_{\lambda^*}$ with $s(m,u\lambda^*)\ne0$.  Here we use
the notation $\lambda^*=-w_0\lambda$ and the fact that
$V_\lambda^*\cong V_{\lambda^*}$.

\begin{definition}\label{definition;brionlagrange}
The \emph{highest-weight set} of $M^\tau$ is the set $\ca C(M^\tau)$
consisting of all quotients $\lambda/n$, with
$\lambda\in\Lambda^*_+\cap\a^*$ and $n>0$, such that there exist $m\in
M^\tau$, $k\in K$, and a $U$-invariant global holomorphic section $s$
of $L^n\boxtimes L_{\lambda^*}$ with $s(m,k\lambda^*)\ne0$.
\end{definition}

\begin{lemma}\label{lemma;weight}
If $M^\tau$ is nonempty\upn, then $\ca C(M^\tau)=\ca C(M)\cap\a^*$.
\end{lemma}

\begin{proof}
The inclusion $\ca C(M^\tau)\subset\ca C(M)\cap\a^*$ is obvious.  Now
let $\xi\in\ca C(M)\cap\a^*$, write $\xi=\lambda/n$ with
$\lambda\in\Lambda^*_+\cap\a^*$ and $n>0$, and choose a $U$-invariant
$s\in\Gamma(M\times U\lambda^*,L^n\boxtimes L_{\lambda^*})$ such that
$s(m,u\lambda^*)\ne0$ for some $m\in M$, $u\in U$.  We need to show
that $m$ can be chosen in $M^\tau$ and $u$ in $K$.

Let us first consider the case $\xi=0$, where we need merely show that
$s(m)\ne0$ for some $m\in M^\tau$.  Here the desired result follows
from Proposition \ref{proposition;lagrangian}, which says that
$M^\tau$ is a totally real Lagrangian: if $s$ vanished on $M^\tau$,
then it would vanish on all of $M$ by the identity principle.

The general case is reduced to the special case by means of the
``shifting trick'': one replaces $M$ with the product $M\times
U\lambda^*$ and $L$ with $L^n\boxtimes L_{\lambda^*}$.  Lemma
\ref{lemma;fix} below implies that $\lambda^*\in\a^*$, so according to
Example \ref{example;orbit} the orbit $U\lambda^*$ has a well-defined
antisymplectic involution, which we denote here by $\tau_{\lambda^*}$,
and $(U\lambda^*)^{\tau_{\lambda^*}} =K\lambda^*\ne\emptyset$ by
\eqref{equation;connect}.  Applying the result of the previous
paragraph, we find a point in
$M^\tau\times(U\lambda^*)^{\tau_{\lambda^*}} =M^\tau\times K\lambda^*$
where $s$ does not vanish, i.e.\ $m\in M^\tau$ and $k\in K$ such that
$s(m,k\lambda^*)\ne0$.
\end{proof}

\begin{lemma}\label{lemma;fix}
Let $w_0^\a$ be the longest element in $W^\a$.  Then
$w_0^\a\lambda=w_0\lambda$ for all $\lambda\in\a^*$.  It follows that
$w_0\in W^\sigma$.
\end{lemma}

\begin{proof}
Let $\lambda\in\a^*_+$.  By definition $w_0^\a$ satisfies
$w_0^\a\lambda\in-\a^*_+\subset-\t^*_+$.  Furthermore,
$w_0^\a\lambda=w\lambda$ for some $w\in W^\sigma$ by Lemma
\ref{lemma;weyl}\eqref{item;normal}, so $w\lambda\in-\t^*_+$.  At the
same time $w_0\lambda\in-\t^*_+$ and therefore $w\lambda=w_0\lambda$,
because $W\lambda$ intersects $-\t^*_+$ in a unique point.  We
conclude that $w_0^\a\lambda=w_0\lambda$.  By Lemma
\ref{lemma;weyl}\eqref{item;normal} this implies $w_0\in\ca
N_W(\a)=W^\sigma$.
\end{proof}

Let $\ca X(M)$ be the Lie algebra of vector fields on $M$.  Since $\ca
X(M)$ has a complex structure, the homomorphism $\ca A\colon\u\to\ca
X(M)$ which maps $\xi$ to $\xi_M$ extends naturally to a homomorphism
of complex Lie algebras,
$$
\ca A^\C\colon\u^\C\map\ca X(M).
$$
Recall that the group of biholomorphic transformations of a compact
complex manifold is a complex Lie group.  It follows from this that
$\ca A^\C$ integrates to an action of the complex reductive group
$U^\C$.  Let $\sigma^\C$ denote the complex linear extension of
$\sigma$ to $\u^\C$ (and also the holomorphic extension of $\sigma$ to
$U^\C$), and let $\gamma_\u$ be complex conjugation on $\u^\C$
relative to $\u$.  Put $\gamma_\g =\sigma^\C\gamma_\u$.  Then
$\gamma_\g$ is a conjugate linear involution of $\u^\C$ because
$\sigma^\C\gamma_\u =\gamma_\u\sigma^\C$.  It is nothing other than
complex conjugation on $\u^\C$ relative to the real form
\begin{equation}\label{equation;dualdecompose}
\g=\k\oplus\p,
\end{equation}
where $\p=i\q$.  It follows from \eqref{equation;semi} that $\ca
A\sigma=\tau\ca A$ and therefore
\begin{equation}\label{equation;complex}
\ca A^\C\gamma_\g=\tau\ca A^\C,
\end{equation}
because $\tau$ is conjugate linear.  Let us write $G$ for the
connected real reductive subgroup of $U^\C$ generated by $\exp\g$ and
$\theta$ for the restriction of $\sigma^\C$ to $G$.  Then $K$ is a
maximal compact subgroup of $G$, $\theta$ is the Cartan involution,
and the symmetric pair $(\g,\k)$ is known as the \emph{dual} of
$(\u,\k)$; cf.\ \cite[Ch.\
V]{helgason;differential-geometry-lie-groups}.  The following result
is analogous to Proposition \ref{proposition;lagrangian}.

\begin{proposition}\label{proposition;orbit}
The involution $\tau$ is $G$-equivariant.  Hence every connected
component of $M^\tau$ is $G$-stable.  For every $m\in M^\tau$ the
intersection $U^\C m\cap M^\tau$ has a finite number of connected
components\upn, each of which consists of a single $G$-orbit.
\end{proposition}

\begin{proof}
It follows from \eqref{equation;complex} that $\tau_*\xi_M=\xi_M$ for
all $\xi\in\g$.  Consequently $\tau(gm)=g\tau(m)$ for all $g\in G$ and
$m\in M$ because $G$ is connected.  This implies that $G$ preserves
every component of $M^\tau$.  That $U^\C m\cap M^\tau$ has finitely
many connected components follows from the fact that it is a real
algebraic variety.  To prove that each component is a single $G$-orbit
it suffices to show that the tangent space to $Gm$ at $m$ is equal to
$T_m(U^\C m)\cap T_mM^\tau$.  A general tangent vector to $U^\C m$ at
$m$ is of the form $\zeta_M=\zeta^+_M+\zeta^-_M$, where
$\zeta^\pm=\xi^\pm+i\eta^\pm$ with $\xi^+$, $\eta^+\in\k$ and $\xi^-$,
$\eta^-\in\q$.  According to \eqref{equation;complex},
$$
\tau_*(\zeta_M)=\gamma_\g(\zeta)_M=(\xi^+-i\eta^+-\xi^-+i\eta^-)_M,
$$
so $\zeta_M$ is tangent to $M^\tau$ if and only if
$\tau_*(\zeta_M)=\zeta_M$ if and only if $(\xi^-+i\eta^+)_M=0$, i.e.\
$\zeta_M=(\xi^++i\eta^-)_M\in T_m(Gm)$.
\end{proof}

\begin{example}[real flag varieties]\label{example;realflag}
Let $M$ be a symmetric coadjoint orbit, $M=U\lambda$ with
$\lambda\in\q^*$.  We saw in Example \ref{example;orbit} that $M^\tau$
consists of the single $K$-orbit $K\lambda$ and in the proof of
Theorem \ref{theorem;orbitquant} that $M$ has a natural compatible
complex structure.  Therefore $G$ acts in a natural way on $M$ and
$M^\tau$ consists of a single $G$-orbit.  As explained in
\cite[\S5]{duistermaat;convexity}, the stabilizer $G_\lambda$ is the
real parabolic subgroup $K_\lambda \exp i\a\exp\lie n$, where $\lie n$
is the nilpotent algebra $\g\cap\bigoplus_{\alpha\in R_-\setminus
R'_-}\C E_\alpha$.  Thus $M^\tau$ is a generalized flag variety of the
real reductive group $G$.  See further Example \ref{example;cross} and
Section \ref{section;application}.
\end{example}

We now explain the relationship between highest-weight sets and
polytopes.

\begin{theorem}\label{theorem;kahlerconvex}
Assume that $M^\tau$ is nonempty.  Then the highest-weight set $\ca
C(M)$ is equal to the set of rational points in $\Delta(M)$\upn, and
$\Delta(M)$ is equal to the closure of $\ca C(M)$ in $\t^*_+$.  The
same statement holds with $M$ replaced by $M^\tau$.  Therefore
$\Delta(M^\tau) =\Delta(M)\cap\a^*$ and $\Delta(M^\tau)$ is a rational
convex polytope.
\end{theorem}

\begin{proof}
The first assertion is proved in
\cite{guillemin-sternberg;convexity;;1982,ness;stratification}; cf.\
also \cite{brion;image}.  The proof of the second assertion follows
the same pattern.  We begin by proving that $\ca C(M^\tau)$ is equal
to the set of rational points in $\Delta(M^\tau)$.  As a special case
we show first that
\begin{equation}\label{equation;zero}
0\in\ca C(M^\tau)\quad\iff\quad0\in\Delta(M^\tau).
\end{equation}
Assume that $0\in\ca C(M^\tau)$.  Pick $n>0$ and a $U$-invariant
$s\in\Gamma(M,L^n)$ such that $s(m)\ne0$ for some $m\in M^\tau$.
Choose a $U$-invariant inner product on $\u$.  Then
\begin{equation}\label{equation;yangmills}
\grad\lVert\Phi(m)\rVert^2=2J\Phi(m)^\flat_{M,m},
\end{equation}
where $\flat\colon\u^*\to\u$ is the linear isomorphism defined by the
inner product, and $\Phi(m)^\flat_M$ is the vector field on $M$
induced by $\Phi(m)^\flat$.  (See e.g.\
\cite[\S6]{kirwan;cohomology-quotients-symplectic}.)  Let $\ca F$ be
the flow of the vector field $-\grad\lVert\Phi\rVert^2$ and let $f$ be
the function $\langle s,s\rangle$.  An easy calculation using
\eqref{equation;yangmills} shows that
\begin{equation}\label{equation;maximum}
\frac{d}{dt}f\bigl(\ca F(t,m)\bigr)=8\pi\bigl\lVert\Phi\bigl(\ca
F(t,m)\bigr)\bigr\rVert^2f\bigl(\ca F(t,m)\bigr)\geq0.
\end{equation}
Let $m_\infty=\lim_{t\to\infty}\ca F_t(m)$.  Then
\eqref{equation;maximum} and the fact that $f(m)>0$ imply that
$f(m_\infty)>0$.  Moreover, $\lim_{t\to\infty}df\bigl(\ca
F(t,m)\bigr)\big/dt=0$ and hence
\begin{equation}\label{equation;zerolevel}
\Phi(m_\infty)=0.
\end{equation}
Furthermore, $\Phi(m)\in\q^*$ because $m\in M^\tau$, so
$i\Phi(m)^\flat\in i\q=\p$.  Together with \eqref{equation;yangmills}
this implies that $\grad\lVert\Phi(m)\rVert^2$ is tangent to the orbit
$Gm$.  This orbit is contained in $M^\tau$ by Proposition
\ref{proposition;orbit}.  The conclusion is that the flow $\ca F$
preserves the submanifold $M^\tau$ and also the $G$-orbits inside it.
This means that $m_\infty\in\overline{Gm}\subset M^\tau$ and hence, by
\eqref{equation;zerolevel}, $0\in\Delta(M^\tau)$.

Conversely, if $0\in\Delta(M^\tau)$ then clearly $0\in\Delta(M)$, so
$0\in\ca C(M)$ by the first part of this theorem.  Hence $0\in\ca
C(M^\tau)$ by Lemma \ref{lemma;weight}.  This concludes the proof of
\eqref{equation;zero}.

To prove that $\lambda/n\in\ca C(M^\tau)$ if and only if
$\lambda/n\in\Delta(M^\tau)$, where $\lambda\in\Lambda^*_+$ and $n>0$,
we use the shifting trick as in the proof of Lemma \ref{lemma;weight}.
We note that $\lambda/n\in\Delta(M^\tau)$ if and only if
$0\in\Delta(M^\tau\times K\lambda^*)$, and that $\lambda/n\in\ca
C(M^\tau)$ if and only if $0\in\ca C(M^\tau\times K\lambda^*)$.  (Here
we have multiplied the symplectic form on $M$ by $n$ and raised its
prequantum line bundle to the $n$-th power).  Now apply the previous
argument to the product $M^\tau\times K\lambda^*$.

Next we show that $\ca C(M^\tau)$ is dense in $\Delta(M^\tau)$.  This
follows from
\begin{equation}\label{equation;dense}
\overline{\ca C(M^\tau)} =\overline{\ca C(M)}\cap\a^*
=\Delta(M)\cap\a^* \supset\Delta(M^\tau),
\end{equation}
where in the first equality we have used that the subspace $\a^*$ of
$\t^*$ is defined over the rationals and in the second equality that
the closure of $\ca C(M)$ is $\Delta(M)$.  This proves the second
statement of the theorem.

The last statement follows immediately from the first two and Lemma
\ref{lemma;weight}.
\end{proof}

\begin{remark}
Let $\sigma_+$ be the involution of $\t^*$ defined by $\sigma_+
=-w'_0\circ\sigma$.  It is easy to see that $\sigma_+$ maps
$\Delta(M)$ into itself, and it follows from Lemma
\ref{lemma;centrallong}\eqref{item;fix} that $\Delta(M^\tau)$ is
contained in $\Delta(M)^{\sigma_+}$.  However, if $w'_0\ne1$, then
$\Delta(M)^{\sigma_+}$ may be strictly larger than $\Delta(M^\tau)$.
\end{remark}

Theorem \ref{theorem;kahlerconvex} can be generalized to certain
singular subsets of $M^\tau$.

\begin{definition}\label{definition;pair}
A \emph{$(U,\sigma)$-pair} in $M$ is a pair $(X,Y)$ of subsets
$X\subset M$ and $Y\subset M^\tau$ subject to the following
requirements: $X$ is a $U$-stable and $\tau$-stable irreducible closed
complex algebraic subvariety of $M$ such that $X\reg$ intersects
$M^\tau$, and $Y$ is the closure of $S$, where $S$ is any connected
component of $X\reg\cap M^\tau$.
\end{definition}

Here $X\reg$ denotes the set of smooth points of $X$.  Standard
results in real algebraic geometry (see e.g.\ \cite[Ch.\
2]{benedetti-risler;real-algebraic-semi}) say that $X\reg\cap M^\tau$
has only finitely many connected components and that $Y$ is a
semi-algebraic subset of $M^\tau$.  It follows from Proposition
\ref{proposition;orbit} and the fact that $X\reg$ is $U^\C$-stable
that $S$, and hence $Y$, are stable under the action of $G$.  It is
obvious how to extend the definitions of $\Delta$ and $\ca C$ to this
setting: one puts $\Delta(X)=\Phi(X)\cap\t^*_+$ and $\Delta(Y)
=\Phi(Y)\cap\t^*_+$; and in the definition of $\ca C(M)$ and $\ca
C(M^\tau)$ one replaces $M$ by $X$, $M^\tau$ by $Y$, and $L$ by
$L|_X$.

\begin{theorem}\label{theorem;singular}
Let $(X,Y)$ be any $(U,\sigma)$-pair in $M$.  Then Theorem
{\rm\ref{theorem;kahlerconvex}} remains valid if we replace $M$ by $X$
and $M^\tau$ by $Y$.
\end{theorem}

\begin{proof}
This is proved in exactly the same way as Theorem
\ref{theorem;kahlerconvex}.  The hypotheses on $X$ and $Y$ imply that
$Y$ contains a Lagrangian submanifold of $X\reg$.  Therefore $Y$ is
Zariski dense in $X$, that is to say, the smallest complex subvariety
of $M$ containing $Y$ is equal to $X$.  This ensures that a section of
$L|_X$ which vanishes on $Y$ vanishes on $X$.  The fact that $Y$ is
$G$-stable implies that the flow $\ca F$ maps $Y$ into itself.  Since
$Y$ is closed the limit of every trajectory in $Y$ is contained in
$Y$.  These facts guarantee that the arguments used to prove Lemma
\ref{lemma;weight} and Theorem \ref{theorem;kahlerconvex} still work
in the present situation.
\end{proof}

For instance, let $X=M$ and let $Y$ be any connected component of
$M^\tau$.  Then $(X,Y)$ is a $(U,\sigma)$-pair.

\begin{corollary}\label{corollary;component}
If $Y$ is a connected component of $M^\tau$\upn, then
$\Delta(Y)=\Delta(M)\cap\a^*$.
\end{corollary}

More interesting examples are obtained by considering any point $m$ in
$M^\tau$ and letting $X$ be the closure of the orbit $U^\C m$.  Then
$X\reg$ contains $U^\C m$, so we can take $Y$ to be the closure of any
component of $U^\C m\cap M^\tau$, for instance $Y=\overline{Gm}$.  The
conclusion is as follows.

\begin{corollary}[orbit closures]\label{corollary;orbitclosure}
$\Delta(\overline{Gm})=\Delta(\overline{U^\C m})\cap\a^*$ for every
$m\in M^\tau$.
\end{corollary}

This is a generalization to the nonabelian case of Theorem 2 in
Atiyah's paper \cite{atiyah;convexity-commuting}.  Brion
\cite{brion;image} proved that for ``most'' $m$ the moment polytope of
the orbit closure $\overline{U^\C m}$ is equal to the full moment
polytope $\Delta(M)$.  We show now that a similar result holds for the
real orbits $Gm$.  This was first pointed out by Duistermaat
\cite{duistermaat;convexity} in the abelian case.  Our proof is
patterned on Brion's.  If $\bb K$ denotes the field $\R$ or $\C$, let
us refer to open subsets in the Zariski topology over $\bb K$ as
\emph{$\bb K$-Zariski open}.

\begin{corollary}\label{corollary;dual}
For all $(U,\sigma)$-pairs $(X,Y)$ there exists a nonempty
$\C$-Zariski open subset $O$ of $X$ such that
$\Delta(Y)=\Delta(\overline{Gm})$ for all $m\in O\cap Y$.  The set
$O\cap Y$ is dense in $Y$.
\end{corollary}

The second statement does not follow automatically from the first,
because $\R$-Zariski open subsets of an irreducible real algebraic
variety are not necessarily dense.  They are always dense, however, if
the variety is nonsingular.

\begin{proof}
Let $\xi_1$, $\xi_2,\dots$, $\xi_k\in\a^*$ be the vertices of the
polytope $\Delta(Y)$.  For each $i$ select
$\lambda_i\in\Lambda^*_+\cap\a^*$ and $n_i>0$ such that
$\xi_i=\lambda_i/n_i$ and pick a nonzero $U$-invariant
$s_i\in\Gamma(X,L^n)$.  Let $O$ be the $\C$-Zariski open subset of $X$
consisting of all $m\in X$ such that $s_i(m)\ne0$ for all $i$.  For
every $m\in O\cap Y$, $s_i(m)\ne0$ and therefore $\xi_i\in\ca
C(\overline{Gm})$ for all $i$.  Using Theorem \ref{theorem;singular}
we conclude that $\Delta(Y)\subset\Delta(\overline{Gm})$.  On the
other hand, it is plain that $\Delta(\overline{Gm})\subset\Delta(Y)$
and so the two polytopes are equal.

The subset $O$ intersects the real algebraic variety $X\reg\cap
M^\tau$ in an $\R$-Zariski open subset and it intersects each
component $S$ of $X\reg\cap M^\tau$ nontrivially since $S$ is a
totally real Lagrangian in $X\reg$.  Added to the fact that $X\reg\cap
M^\tau$ is nonsingular, this implies that $O\cap X\reg\cap M^\tau$ is
dense in $X\reg\cap M^\tau$.  Consequently $O\cap Y$ is dense in $Y$.
\end{proof}

\begin{example}[products of orbits]\label{example;cross}
Let $M=U\lambda\times U\mu$ be the product of the coadjoint orbits
through $\lambda$, $\mu\in\a^*_+$.  Then $M^\tau=K\lambda\times K\mu$.
Let us assume for simplicity that $\lambda$ and $\mu$ are generic, so
that $U_\lambda=U_\mu=\ca Z_K(A)=K'$.  According to Example
\ref{example;realflag} we can identify $M^\tau$ with $G/B\times G/B$,
where $B$ is the real Borel subgroup $K'\exp i\a\exp\lie n$ of $G$.
We see immediately that the $G$-orbits on $G/B\times G/B$ are in a
natural one-to-one correspondence with the $B$-orbits on $G/B$.
Moreover, this correspondence preserves the closure relation among
orbits.  Applying the Bruhat decomposition (Theorem 1.4 in \cite[Ch.\
IX]{helgason;differential-geometry-lie-groups}) we obtain the
following result: the $G$-orbits $\ca O_w$ on $M^\tau$ are of the form
$$
\ca O_w=G(\bar 1,\bar w)
$$
for $w\in W^\a$; and $\ca O_v$ is in the closure of $\ca O_w$ if and
only if $v\le w$.  Here $\bar w$ denotes the coset $n\bmod B\in G/B$
determined by any representative $n\in\ca N_K(A)$ of $w\in W^\a$, and
``$\le$'' denotes the Bruhat-Chevalley ordering on $W^\a$.  Thus there
is exactly one closed orbit, the diagonal $\ca O_1$, and exactly one
open orbit, $\ca O_{w_0^\a}$.  It is clear that $\Delta(\ca O_1)
=\{\lambda+\mu\}$ and $\Delta(\bar\ca O_{w_0^\a})=\Delta(M^\tau)$.
See Section \ref{section;application} for a complete calculation of
$\Delta(M^\tau)$.  Let us here calculate it and its subpolytopes
$\Delta(\bar\ca O_w)$ under the simplifying assumption that
$\Delta(M^\tau)$ is entirely contained in the interior of $\t^*_+$.
Then Theorem \ref{theorem;convex} implies that all vertices of
$\Delta(M^\tau)$ are images under $\Phi$ of $A$-fixed points on
$K\lambda\times K\mu$.  According to Lemma
\ref{lemma;weyl}\eqref{item;suborbit} the set of $A$-fixed points is
equal to $W^\a\lambda\times W^\a\mu$, and therefore
$$
\Delta(M^\tau)=\hull\{\,u\lambda+v\mu\mid\text{$u$, $v\in W^\a$,
$u\lambda+v\mu$ dominant}\,\}.
$$
Likewise, the vertices of $\Delta(\bar\ca O_w)$ are images of
$A$-fixed points on $\bar\ca O_w=\bigcup_{v\le w}\ca O_v$.  An easy
computation yields that the $A$-fixed points on $\ca O_v$ are the
points $n(\bar 1,\bar v)$ with $n\in\ca N_K(A)$.  Thus
$$
\Delta(\bar\ca O_w)=\hull\{\,u(\lambda+v\mu)\mid\text{$u$, $v\in
W^\a$, $v\le w$, $u(\lambda+v\mu)$ dominant}\,\}.
$$
\end{example}

\section{The affine case}\label{section;affine}

Before we treat the general case we prove a version of the Lagrangian
convexity theorem for affine varieties.  This is not really a special
case of Theorem \ref{theorem;convex}, because the moment map on an
affine variety is not proper, but it is needed in the proof of Theorem
\ref{theorem;convex}.  Also, as in the projective case, the result
holds for singular as well as nonsingular varieties, for instance
closures of orbits of the noncompact real reductive group $G$.  The
idea of the proof is the same as in the projective case.

Let $(V,\tau)$ be a (finite-dimensional) unitary $(U,\sigma)$-module
and let $M$ be a $U$-stable and $\tau$-stable irreducible affine
subvariety of $V$.

\begin{definition}\label{definition;affinelagrange}
The \emph{highest-weight set} of $M$ is the set $\ca C(M)$ consisting
of all dominant weights $\lambda$ such that the irreducible
representation $V_\lambda$ occurs in the coordinate ring $\C[M]$.  The
\emph{highest-weight set} of $M^\tau$ is the set $\ca C(M^\tau)$ of
all $\lambda\in\Lambda^*_+\cap\a^*$ such that there exist $m\in M$,
$k\in K$, and a $U$-invariant global algebraic section $s$ of
$L\boxtimes L_{\lambda^*}$ with $s(m,k\lambda^*)\ne0$.
\end{definition}

Here $L$ denotes the trivial line bundle on $M$ with the trivial lift
of the $U$-action.  The proof of the following result is similar to
that of Lemma \ref{lemma;weight}.

\begin{lemma}\label{lemma;affineweight}
If $M^\tau$ contains a smooth point of $M$\upn, then $\ca
C(M^\tau)=\ca C(M)\cap\a^*$.
\qed
\end{lemma}

Proposition \ref{proposition;orbit} is valid exactly as stated in the
present situation.  Here is the analogue of Theorem
\ref{theorem;kahlerconvex}.

\begin{theorem}\label{theorem;affineconvex}
Assume that $M^\tau$ contains a smooth point of $M$.  Then
$\hull_\Q\ca C(M)$ is equal to the set of rational points in
$\Delta(M)$\upn, and $\Delta(M)$ is equal to the closure of
$\hull_\Q\ca C(M)$ in $\t^*_+$.  The same statement holds with $M$
replaced by $M^\tau$.  Therefore $\Delta(M^\tau) =\Delta(M)\cap\a^*$
and $\Delta(M^\tau)$ is a rational convex polyhedral cone.
\end{theorem}

Here ``$\hull_\Q$'' denotes the convex hull of a subset of a vector
space over $\Q$.  The assumption that $M$ is irreducible implies that
$\ca C(M)$ is closed under addition, and therefore $\hull_\Q\ca C(M)$,
resp.\ $\hull_\Q\ca C(M^\tau)$, is the set of all positive rational
multiples of elements of $\ca C(M)$, resp.\ $\ca C(M^\tau)$.

\begin{proof}
See \cite[Theorem 4.9]{sjamaar;convexity} for the first statement.
Next, let us show that $\hull_\Q\ca C(M^\tau)$ is a subset of
$\Delta(M^\tau)$.  Let $\xi\in\hull_\Q\ca C(M^\tau)$; then
$\lambda=n\xi\in\ca C(M^\tau)$ for some positive integer $n$.
Consider the product $M'=M\times U\lambda^*$ equipped with the
symplectic form $\omega'=n\omega+\omega_{\lambda^*}$ and the moment
map $\Phi'=n\Phi+\iota_{\lambda^*}$.  Here $\omega_\lambda$ denotes
the symplectic form on $U\lambda$ and $\iota_\lambda$ the inclusion of
$U\lambda$ into $\u^*$.  We need to show that
$0\in\Delta\bigl((M')^\tau\bigr)$.  This is established by picking a
$U$-invariant global algebraic section $s$ of $L'=L\boxtimes
L_{\lambda^*}$ and a point $m'=(m,k\lambda)\in(M')^\tau$ with
$s(m')\ne0$, and considering the trajectory of $m'$ under the flow
$\ca F'$ of $-\grad\lVert\Phi'\rVert^2$.  It is proved in
\cite{sjamaar;convexity} that $m'_\infty=\lim_{t\to\infty}\ca
F'_t(m')$ exists.  One then shows as in the proof of Theorem
\ref{theorem;kahlerconvex} that $m'_\infty\in(M')^\tau$ and
$\Phi'(m'_\infty)=0$.  (One uses the Hermitian metric on $L$ defined
by the Gaussian $\exp(-\pi n\lVert\cdot\rVert^2)$.)

Now assume that $\xi$ is a rational point in $\Delta(M^\tau)$.  Then
$\xi\in\Delta(M)\cap\a^*$, so, by the first part of the theorem
$\xi\in\hull_\Q\ca C(M)\cap\a^*$, which is equal to $\hull_\Q\ca
C(M^\tau)$ by Lemma \ref{lemma;affineweight}.  We have shown that
$\Delta(M)\cap(\Lambda^*\otimes\Q)=\hull_\Q\ca C(M^\tau)$.

The fact that $\hull_\Q\ca C(M^\tau)$ is dense in $\Delta(M^\tau)$ is
proved as in \eqref{equation;dense}.

The last assertion of the theorem follows from the first two plus
Lemma \ref{lemma;affineweight}.
\end{proof}

\begin{corollary}
The moment cones $\Delta(M)$ and $\Delta(M^\tau)$ are independent of
the $U\rtimes\{\pm1\}$-equivariant affine embedding of $M$ into the
$(U,\sigma)$-module $V$.
\qed
\end{corollary}

We leave it to the reader to state and prove analogues of Corollaries
\ref{corollary;orbitclosure} and \ref{corollary;dual}, and finish this
section with a necessary criterion for the origin to be an extreme
point of the moment cone.

\begin{lemma}\label{lemma;extreme}
Assume that $\Delta(M^\tau)$ is a proper cone.  Then
$A\subset[U,U]\,U_m$ for all $m$ such that $U^\C m$ is closed in $M$.
\end{lemma}

\begin{proof}
Consider the homogeneous variety $X=U^\C/[U^\C,U^\C]\,U^\C_m$.  The
closedness of the orbit $U^\C m$ ensures that $(U^\C)_m$ is the
complexification of $U_m$ and that $X$ is affine.  According to the
proof of \cite[Theorem 4.25]{sjamaar;convexity}, $\Delta(X)$ is
contained in $\Delta(M)$ and is equal to the kernel of the projection
$p\colon\u^*\to([\u,\u]+\u_m)^*$.  Therefore, if
$\Delta(M^\tau)=\Delta(M)\cap\a^*$ is a proper cone, then $\ker
p\cap\a^*=\{0\}$.  Hence $p$ induces an injection
$\a^*\hookrightarrow([\u,\u]+\u_m)^*$, and so $\a\subset[\u,\u]+\u_m$.
This is equivalent to $A\subset[U,U]\,U_m$.
\end{proof}

\section{Normal forms}\label{section;normal}

In this section we extend to the category of Hamiltonian
$(U,\sigma)$-manifolds a well-known local normal form theorem of
Guillemin and Sternberg \cite{guillemin-sternberg;normal} and Marle
\cite{marle;modele-action}.  Because this is mostly a routine
exercise, which amounts to inserting $\sigma$s and $\tau$s in the
right places, we omit the proof.  The techniques involved are a
straightforward combination of those used in \emph{loc.\ cit.}, Meyer
\cite{meyer;hamiltonian-discrete}, and Duistermaat
\cite{duistermaat;convexity}.

Let $H$ be a closed (but not necessarily connected) $\sigma$-stable
subgroup of $U$ and let $V$ be a finite-dimensional unitary
$(H,\sigma)$-module, considered as a Hamiltonian $(H,\sigma)$-manifold
with momentum map $\Phi_V$ as in Example \ref{example;linear}.
According to Examples \ref{example;cotangent} and
\ref{example;restrictionproduct} the product $T^*U\times V$ is a
Hamiltonian $(U\times H,\sigma\times\sigma)$-manifold, where $U$ acts
on $T^*U$ by left multiplication, $H$ acts on $T^*U$ by right
multiplication and linearly on $V$.  Therefore, by Example
\ref{example;quotient}, the symplectic quotient $(T^*U\times V)\qu H$
is a Hamiltonian $(U,\sigma)$-manifold.  There is a naturally defined
$U$-equivariant diffeomorphism between $(T^*U\times V)\qu H$ and the
associated bundle
$$
F(H,V)=U\times^H(\h^0\times V),
$$
which depends only on the choice of a connection on the principal
$H$-bundle $U\to U/H$.  Here $\h^0$ denotes the annihilator of $\h$ in
$\u^*$.  (See e.g.\ \cite[Appendix A]{meinrenken-sjamaar;singular}.)

Now let $M$ be an arbitrary Hamiltonian $(U,\sigma)$-manifold and let
$m\in\Phi\inv(0)$.  Let $H$ be the stabilizer of $m$ and $V
=T_m(Um)^\omega/T_m(Um)$ the symplectic slice at $m$.  It is expedient
to choose a $U$-invariant and $\sigma$-anti-invariant compatible
almost complex structure on $M$; then the symplectic slice is a
unitary representation in a natural way.  Moreover, if $\tau(m)=m$,
the subgroup $H$ is $\sigma$-stable, the subspace $\h^0$ of $\u^*$ is
$\sigma$-stable, and $V$ is a unitary $(H,\sigma)$-module.  We can now
state the local normal form theorem.

\begin{theorem}[symplectic slices]\label{theorem;slice} 
Let $m\in\Phi\inv(0)\cap M^\tau$ and $H=U_m$.  Let $V$ be the
symplectic slice at $m$.  Then there exists a
$U\rtimes\{\pm1\}$-stable open neighbourhood of $m$ in $M$ which is
isomorphic as a Hamiltonian $(U,\sigma)$-manifold to a neighbourhood
of the zero section in $F(H,V)$.  The isomorphism maps $m$ to
$[1,0,0]$\upn; the moment map and the involution on $F(H,V)$ are given
by
\begin{align*}
\Phi\colon[u,\xi,v] &\longmapsto u\bigl(\xi+\Phi_V(v)\bigr),\\
\tau\colon[u,\xi,v] &\longmapsto[\sigma(u),-\sigma(\xi),\tau_V(v)],
\end{align*}
respectively\upn, where $\tau_V$ denotes the involution induced on
$V$.
\qed
\end{theorem}

One interesting consequence of this result concerns the symplectic
quotient $M\qu U$.  The involution $\tau$ preserves the fibre
$\Phi\inv(0)$ and maps $U$-orbits to $U$-orbits.  It therefore
descends to an involution on $M\qu U$, which we denote also by $\tau$.
If $U$ acts freely on $\Phi\inv(0)$, then $M\qu U$ is a symplectic
manifold and the involution is antisymplectic.  Therefore the
fixed-point set $(M\qu U)^\tau$ is a Lagrangian submanifold.  The
\emph{Lagrangian quotient} or \emph{Lagrangian reduced space}
$M^\tau\qu U$ is by definition the space obtained by projecting
$UM^\tau\cap\Phi\inv(0)$ into $M\qu U$.

\begin{corollary}\label{corollary;quotient}
If $U$ acts freely on $\Phi\inv(0)$\upn, then the Lagrangian quotient
$M^\tau\qu U$ is a union of connected components of the Lagrangian
$(M\qu U)^\tau$.
\end{corollary}

\begin{proof}
Clearly both $M^\tau\qu U$ and $(M\qu U)^\tau$ are closed and
$M^\tau\qu U\subset(M\qu U)^\tau$, so it suffices to show that
$M^\tau\qu U$ is a Lagrangian submanifold.  Let $m\in M^\tau$.  We
compute in the model $F(H,V)$.  Writing $\h^0=\lie m$ we have an
inclusion
\begin{equation}\label{equation;locallagrangian}
U^\sigma\times^{H^\sigma}(\lie m^{-\sigma}\times V^\tau)\subset
F(H,V)^\tau
\end{equation}  
of closed Lagrangian submanifolds of $F(H,V)$.  Upon taking quotients
we obtain
\begin{equation}\label{equation;quotientlagrangian}
\bigl(V^\tau\cap\Phi_V\inv(0)\bigr)\big/H^\sigma\subset F(H,V)^\tau\qu
U.
\end{equation}
Now we use the assumption that $U$ acts freely.  This implies that
$H=\{1\}$, $\Phi_V=0$, and also that the inclusions
\eqref{equation;locallagrangian} and
\eqref{equation;quotientlagrangian} are equalities.  Therefore
$F(H,V)\qu U=V$ and $F(H,V)^\tau\qu U =V^\tau$.  The upshot is that
$F(H,V)^\tau\qu U$ is a Lagrangian submanifold of $F(H,V)\qu U$.
\end{proof}

\begin{remark}[singular quotients]\label{remark;singular}
If we drop the assumption that $U$ acts freely on $\Phi\inv(0)$, the
symplectic quotient $M\qu U$ is a symplectic stratified space.  An
argument analogous to the proof of Corollary \ref{corollary;quotient}
shows that the decomposition of $M\qu U$ into symplectic strata
induces a decomposition of the Lagrangian quotient $M^\tau\qu U$ into
smooth manifolds, each of which is Lagrangian in the ambient
symplectic stratum of $M\qu U$.  Moreover, $M^\tau\qu U$ always
intersects the open stratum of $M\qu U$, so that $\dim M^\tau\qu
U=\frac1{2}\dim M\qu U$.
\end{remark}

Another consequence of the normal form theorem is the fact that $M$ is
locally, near the orbit through any point $m\in
M^\tau\cap\Phi\inv(0)$, isomorphic to an affine algebraic Hamiltonian
$(U,\sigma)$-manifold.  As before, let $H=U_m$ and let $V$ be the
symplectic slice at $m$.  Consider the associated bundle
$U^\C\times^{H^\C}V$, which is an affine variety with a natural
algebraic action of $U^\C$.

\begin{lemma}\label{lemma;affinebundle}
There exist a finite-dimensional unitary $(U,\sigma)$-module
$(E,\tau)$ and a $U$-equivariant algebraic embedding $\iota\colon
U^\C\times^{H^\C}V\hookrightarrow E$ such that the image
$X=\iota(U^\C\times^{H^\C}V)$ is $\tau$-stable.  Furthermore\upn, if
$x=\iota([1,0])$\upn, then $\Phi_E(x)=0$ and the symplectic slice at
$x$ \upn(relative to the symplectic structure on $X$ induced by the
constant symplectic form on $E$\upn) is isomorphic to $V$.
\end{lemma}

\begin{proof}[Outline of proof]
Without the $\tau$, this is \cite[Example 5.5]{sjamaar;convexity}.
Let us show briefly how to adapt the argument to the present context.

We use the fact that for every compact Lie group $\ca G$ and closed
subgroup $\ca H$ there exists a finite-dimensional real $\ca G$-module
with a vector in it whose stabilizer is equal to $\ca H$.  Taking $\ca
G=U\rtimes\{\pm1\}$ and $\ca H=H\rtimes\{\pm1\}$, where $\{\pm1\}$
acts on $U$ and $H$ by means of the involution $\sigma$, gives us a
real $U$-module $E_1^\R$ that is equipped with an involution
$\tau_1^\R$ and contains a vector $e_1$ such that $U_{e_1}=H$,
$\tau_1^\R(e_1)=e_1$, and $\tau_1^\R(ue)=\sigma(u)e$ for all $u\in U$
and $e\in E_1^\R$.  Let $E_1$ be the complexification of $E_1^\R$ and
$\tau_1$ the \emph{conjugate} linear extension of $\tau_1^\R$ to
$E_1$, i.e.\ the complex linear extension composed with complex
conjugation.  Then $(E_1,\tau_1)$ is a unitary $(U,\sigma)$-module,
and the map $U^\C\to E_1$ which sends $g$ to $ge_1$ descends to a
closed equivariant algebraic embedding $\iota_1\colon
U^\C/H^\C\hookrightarrow E_1$.  The image of $\iota_1$ is
$\tau_1$-stable: if $g\in U^\C$, then $g=u\exp(i\xi)$ with $u\in U$
and $\xi\in\u$, so
$$
\tau_1(ge_1)=\tau_1\bigl(u\exp(i\xi)e_1\bigr)
=\sigma\bigl(u\exp(-i\xi)\bigr)\tau_1(e_1)\in U^\C e_1.
$$

Moreover, there exists a (complex, finite-dimensional)
$(U,\sigma)$-module $(E_2,\tau_2)$ together with an embedding
$\iota_2\colon V\hookrightarrow E_2$ which intertwines the $H$-actions
and the involutions.  To prove this, by virtue of the structure
theorem, Addendum \ref{addendum;classification}, we may assume $V$ to
be either of type \eqref{equation;irreducible} or
\eqref{equation;dup}.  (See Remark \ref{remark;disconnect} if $H$ is
not connected.)  Let $F$ be a (complex, finite-dimensional)
$(U,\sigma)$-module containing a vector whose stabilizer is $U$.  Let
$\lambda$ be a dominant weight of $U$.  An argument based on
Peter-Weyl (see e.g.\ \cite{janich;differenzierbare;;1968}) shows that
there exists an $H$-equivariant embedding $j\colon
V_\lambda\hookrightarrow E_2$, where $E_2$ is of the form
$F\otimes\dots\otimes F\otimes F^*\otimes\dots\otimes F^*$.  Observe
that $E_2$ is a $(U,\sigma)$-module with involution $\tau_2$ in a
natural way.  Assume that $V=V_\lambda$ is of type
\eqref{equation;irreducible}.  Put $\iota_2=\frac1{2}(j+\tau_2\circ
j\circ\tau)$.  Then $\iota_2$ is $H$-equivariant and intertwines the
involutions.  By Schur's Lemma, $\iota_2$ is either injective or $0$.
In the first case we are done; in the second case we replace $\tau_2$
by $-\tau_2$ and $\iota_2$ by $\frac1{2}(j-\tau_2\circ j\circ\tau)$.
Now assume that $V=V_\lambda\dup$ is of type \eqref{equation;dup}.
Then the intersection of $j(V_\lambda)$ and
$\tau\bigl(j(V_\lambda)\bigr)$ is $\{0\}$, so
$\iota_2=j\oplus\tau_2\circ j\circ\tau$ defines an $H$-equivariant
embedding of $V\cong V_\lambda\oplus V_{\sigma_+(\lambda)}$ into $E_2$
which intertwines $\tau$ and $\tau_2$.

Finally, we embed the associated bundle into $E_1\oplus E_2$ by
putting $\iota([g,v])=\bigl(ge_1,g\iota_2(v)\bigr)$.
\end{proof}

The upshot of Theorem \ref{theorem;slice} and Lemma
\ref{lemma;affinebundle} is as follows.

\begin{corollary}\label{corollary;localaffine}
Under the hypotheses of Theorem {\rm\ref{theorem;slice}} there exists
a $U\rtimes\{\pm1\}$-stable open neighbourhood of $m$ in $M$ which is
isomorphic as a Hamiltonian $(U,\sigma)$-mani\-fold to a neighbourhood
of $x$ in $X$.
\end{corollary}

\section{Proof of the main theorem}\label{section;proof}

Let $M$ be an arbitrary $(U,\sigma)$-manifold and let $m\in M^\tau$.
Assume that $\lambda=\Phi(m)$ is in the fundamental Weyl chamber.  The
positive roots of $(U_\lambda,T)$ are by definition those positive
roots of $(U,T)$ which are perpendicular to $\lambda$, and the
fundamental Weyl chamber $\t^*_{+,\lambda}$ of $U_\lambda$ is chosen
accordingly.  Let $H=U_m\subset U_\lambda$ be the stabilizer of $m$
and
$$
V=T_m(Um)^\omega\big/\bigl(T_m(Um)\cap T_m(Um)^\omega\bigr)
$$
the symplectic slice at $m$.  Consider the affine
$(U_\lambda,\sigma)$-variety $X=U_\lambda^\C\times^{H^\C}V$ and its
moment cone $\Delta(X)\subset\t^*_{+,\lambda}$.

\begin{definition}\label{definition;local}
The \emph{local moment cone} of $M$ at $m$ is
$\Delta_m=\lambda+\Delta(X)$.  The \emph{local moment cone} of
$M^\tau$ at $m$ is $\Delta^\tau_m=\lambda+\Delta(X^\tau)$.
\end{definition}

It follows from Theorem \ref{theorem;affineconvex} that
$\Delta^\tau_m=\Delta_m\cap\a^*$.  Here is a local version of Theorem
\ref{theorem;convex}.

\begin{theorem}\label{theorem;localcone}
For every sufficiently small $U$-stable and $\tau$-stable open
neighbourhood $O$ of $m$ the set $\Delta(O)$ is a neighbourhood of the
vertex $\lambda$ in $\Delta_m$\upn, and $\Delta(O^\tau)$ is a
neighbourhood of the vertex $\lambda$ in $\Delta^\tau_m$.
\end{theorem}

\begin{proof}
The first statement is \cite[Theorem 6.5]{sjamaar;convexity}.  We now
prove the second statement.  Under the assumption that $\lambda=0$
this follows from Corollary \ref{corollary;localaffine}, together with
the following observation: if $O$ is a $U$-stable and $\tau$-stable
neighbourhood of $x$ in $X$, then $\Delta(O^\tau)$ is a neighbourhood
of the vertex in $\Delta(X^\tau)$.  This would be automatic if the
moment map on $X$ were proper, which it usually is not.  However, we
can factor $\Phi$ into a proper moment map and a linear projection.
The $S^1$-action on $X$ defined by complex scalar multiplication on
the slice $V$ commutes with the $U$-action and is Hamiltonian with
moment map $\phi([g,v])=-\frac1{2}\lVert v\rVert^2$.  It is not hard
to show that the joint moment map $\Phi\times\phi\colon
X\to\u^*\times\R$ is proper.  Clearly $\Phi=\pi\circ(\Phi\times\phi)$,
where $\pi\colon\u^*\times\R\to\u^*$ is the projection onto the first
factor.  For $x\in X$ define $\Phi_+(x)$ to be the unique intersection
point of the coadjoint orbit $U\Phi(x)$ with the Weyl chamber
$\t^*_+$.  Being proper, $\Phi_+\times\phi$ maps $O^\tau$ onto a
neighbourhood of the vertex of the $U\times S^1$-moment cone of
$X^\tau$ and, being linear, $\pi$ projects this neighbourhood onto a
neighbourhood of the vertex in $\Delta(O^\tau)$.

The general case can be reduced to the case $\lambda=0$ by means of
symplectic cross-sections.  Let $\s$ be the (relatively open) face of
$\t^*_+$ containing $\lambda$.  Let $\eu S_\s=U_\s\cdot\star\s$, where
$\star\s$ stands for the open star $\bigcup_{\eu r\suc\s}\eu r$ of
$\s$ and $U_\s$ is the centralizer of $\s$.  Then $\eu S_\s$ is a
maximal slice for the coadjoint action at $\lambda$.  Its preimage
$M_\s=\Phi\inv(\eu S_\s)$ is the \emph{symplectic cross-section} at
$\s$.  (See e.g.\ \cite[Section 3.2]{meinrenken-sjamaar;singular}.)
The fact that $\lambda\in\s\cap\a^*$ implies that $\s$, and hence
$\star\s$, are stable under $\sigma_+$.  Moreover, the longest element
$w_0'\in W'$ has a representative $k_0$ which is in $K'\subset U_\s$.
It follows from this that
$-\sigma(u\xi)=\sigma(u)k_0\sigma_+(\xi)\in\eu S_\s$ for all $u\in
U_\s$ and $\xi\in\star\s$, i.e.\ $\eu S_\s$ is stable under $-\sigma$.
Consequently $M_\s$ is $\tau$-stable and is therefore a Hamiltonian
$(U_\s,\sigma)$-manifold in its own right.  Moreover, $\Delta(M_\s)
=\Delta(M)\cap\star\s$ and $\Delta(M_\s^\tau)
=\Delta(M^\tau)\cap\star\s$ by equivariance of the moment map.  It
therefore suffices to prove the result for $M_\s$ instead of $M$.
Since $\lambda$ is annihilated by $U_\s$, we can subtract $\lambda$
from the moment map to get $\Phi(m)=0$.
\end{proof}

\begin{proof}[Proof of Theorem {\rm\ref{theorem;convex}}]
It is clear that $\Delta(M^\tau)$ is contained in $\Delta(M)\cap\a^*$.
Let $\lambda\in\Delta(M^\tau)$.  According to \cite[Theorem
6.5]{sjamaar;convexity} the properness of $\Phi$ implies that the
local moment cone $\Delta_m$ is the same for all
$m\in\Phi\inv(\lambda)$.  Hence, by Theorem
\ref{theorem;affineconvex}, the cone $\Delta^\tau_m$ is the same for
all $m\in\Phi\inv(\lambda)\cap M^\tau$.  Together with Theorem
\ref{theorem;localcone} this implies that a small neighbourhood of
$\lambda$ in $\Delta(M^\tau)$ is equal to a small neighbourhood of
$\lambda$ in $\Delta^\tau_m$.  In particular, $\Delta(M^\tau)$ is
locally convex.  It is also closed because $\Phi$ is proper, and hence
it is globally convex.  A convex set is the intersection of all convex
cones containing it, so
$$
\Delta(M^\tau)=\bigcap_{m\in M^\tau}\Delta^\tau_m =\bigcap_{m\in
M^\tau}\Delta_m\cap\a^*\subset\bigcap_{m\in M}\Delta_m\cap\a^*
=\Delta(M)\cap\a^*,
$$
and hence $\Delta(M^\tau)=\Delta(M)\cap\a^*$.  This proves
\eqref{item;convex}.

\eqref{item;vertex} follows from Lemma \ref{lemma;extreme} and Theorem
\ref{theorem;localcone}.
\end{proof}

\section{Applications}\label{section;application}

\subsection*{Cotangent bundles}

Let $S$ be a connected manifold equipped with a $U$-action and an
involution $\tau$ such that $\tau(us)=\sigma(u)\tau(s)$.  Let $M=T^*S$
be the $(U,\sigma)$-manifold defined in Example
\ref{example;cotangent}.  The moment map on $M$ is not proper unless
$S$ is a homogeneous $U$-manifold.  Nevertheless, Theorem
\ref{theorem;localcone} enables us to describe the images of $M$ and
$M^\tau=T^*_{S^\tau}S$.  Assume that $S^\tau$ is nonempty and take any
point $s$ in it.  Let $H$ be the stabilizer $U_s$ and
$V^\R=T_sS/T_s(Us)$ the slice at $s$.  Then $\tau$ induces an
involution $\tau^\R$ on $V^\R$ and the symplectic slice $V$ at $s$ is
isomorphic to the complexification of $V^\R$, equipped with the
obvious $U$-action and the conjugate linear extension of $\tau^\R$.
The affine local model $U^\C\times^{H^\C}V$ is isomorphic to the
complexification of the real affine $U$-variety $U\times^HV^\R$.

\begin{theorem}\label{theorem;cotangent}
The set $\Delta(T^*S)$ is a rational convex polyhedral cone.  It is
equal to the moment cone of the affine $U$-variety
$U^\C\times^{H^\C}V$.  Moreover\upn, $\Delta(T^*_{S^\tau}S)$ is equal
to $\Delta(T^*S)\cap\a^*=\Delta(U^\C\times^{H^\C}V)\cap\a^*$ and is
therefore likewise a rational convex polyhedral cone.
\end{theorem}

\begin{proof}
See \cite[Theorem 7.6]{sjamaar;convexity} for the first statement.  To
prove the second statement, observe that the moment map is homogeneous
in the cotangent directions and therefore $\Delta(T^*_{S^\tau}S)$ is
equal to the cone spanned by $\Delta(O)$, where $O$ is any open subset
of $T^*_{S^\tau}S$ containing the zero section.  Now apply Theorem
\ref{theorem;localcone}.
\end{proof}

\subsection*{Coadjoint orbits and subgroups}

Let $M$ be a coadjoint orbit of $U$.  In
\cite{berenstein-sjamaar;coadjoint} a description was given of the
moment polytope of $M$ considered as a Hamiltonian manifold for the
action of a \emph{subgroup} of $U$.  In this section we combine the
results of \cite{berenstein-sjamaar;coadjoint} with Theorem
\ref{theorem;convex} to describe the moment polytope of the real flag
variety $M^\tau$ relative to the action of a subgroup.

More specifically, let $\ti U$ be a closed connected subgroup of $U$
and denote the inclusion map $\ti U\to U$ by $f$.  Assume that $\ti U$
is $\sigma$-stable and denote the restriction of $\sigma$ to $\ti U$
by $\ti\sigma$.  Then we have a decomposition $\ti\u=\ti\k\oplus\ti\q$
analogous to \eqref{equation;decompose}.  Let $\ti Q=\{\,\ti u\ti
\sigma(\ti u)\inv\mid\ti u\in\ti U\,\}$ be the symmetric space and let
$\ti A$ be a maximal torus in $\ti Q$ containing $1$.  Let $\ti T$ be
a maximal torus in $\ti U$ containing $A$ and let $T$ be a maximal
torus in $U$ containing $\ti T$.  It is not necessarily possible to
choose $T$ in such a way that $A=T\cap Q$ is a maximal torus in $Q$,
but for simplicity we shall assume this to be the case.  (It is
certainly possible in the examples we discuss below.)  Let
$f^*\colon\t^*\to\ti\t^*$ be the projection induced by $f$.  Then
$f^*$ maps $\a^*$ onto $\ti\a^*$.

We assume the orbit $M$ to be symmetric, i.e.\ $M=K\lambda$ with
$\lambda\in\a^*$, so that $M^\tau=K\lambda$, and we consider $M$ as a
Hamiltonian $(\ti U,\ti\sigma)$-manifold.  The following inequalities
for $\Delta(M^\tau)$ are now an immediate consequence of Theorem
\ref{theorem;convex} and \cite[Theorem
3.2.1]{berenstein-sjamaar;coadjoint}.

\begin{theorem}\label{theorem;orbit}
Let $(\ti\lambda,\lambda)\in\ti\a^*_+\times\a^*_+$.  Then
$\ti\lambda\in\Delta(M^\tau)$ if and only if
\begin{equation}\label{equation;polyinequal}
\ti w\inv\ti\lambda\in f^*(w\inv\lambda-v\ca C)
\end{equation}
for all triples $(\ti w,w,v)\in\ti W\times W\times W\rel$ such that
$\ti\sigma_{\ti w}$ is contained in $\phi^*(v\sigma_{wv})$.
\qed
\end{theorem}

Here $W\rel$ is the \emph{relative Weyl set}, which is defined as
follows.  The \emph{compatible Weyl set} $W\com$ is the set of all
$w\in W$ such that $\dim(w\t_+\cap\ti\t_+)=\dim\ti\t_+$, where $\t_+$
and $\ti\t_+$ are the respective Weyl chambers in $\t$ and $\ti\t$.
Observe that $W\com$ is stable under left multiplication by $\ca
Z_W(\ti\t)$, the centralizer of $\ti\t$ in $W$.  The relative Weyl set
is then the set of shortest representatives of the left $\ca
Z_W(\ti\t)$-cosets in $W\com$.  For each $\ti w$, $w$ and $v$ the
condition \eqref{equation;polyinequal} represents a finite number of
linear inequalities.  Let $X=U/T$ and $\ti X=\ti U/\ti T$ be the flag
varieties of $U$ and $\ti U$, respectively.  Then $\sigma_w\in
H^{2l(w)}(X)$ denotes the Schubert class associated with $w\in W$ and
$\phi\colon\ti X\hookrightarrow X$ the embedding induced by $f$.  We
say that a Schubert class $\ti\sigma_{\ti w}$ is \emph{contained} in a
cohomology class $c$ if it appears in the expansion of $c$ in terms of
the Schubert basis.  See \cite{berenstein-sjamaar;coadjoint} for a
more detailed discussion.

Let $\Delta(T_K^*U)\subset\ti\a^*_+\times\a^*_+$ be the moment cone of
the conormal bundle $T_K^*U$, considered as a $\ti U\times
U$-manifold.  The polytopes $\Delta(M^\tau)$ can be obtained by
slicing this cone horizontally.  The inequalities
\eqref{equation;polyinequal} also describe $\Delta(T_K^*U)$ and thus
Theorem \ref{theorem;orbit} can be regarded as a more explicit version
of Theorem \ref{theorem;cotangent} for $S=U$.

\subsection*{The maximal torus}

The easiest special case of Theorem \ref{theorem;orbit} is that of
$\ti U=T$, the maximal torus of $U$.  Here of course $\ti T=T$, $\ti
W=\{1\}$ and $f^*$ is the identity map.  Furthermore, the flag variety
$\ti X$ is a point, so $\phi^*$ is the trivial homomorphism
$H^\bu(X,\Z)\to\Z$.  It follows that $\phi^*(v\sigma_{wv})\ne0$ if and
only if $v=w\inv$.  Hence for every $\lambda\in\a^*_+$ the points
$\ti\lambda$ in the moment polytope of the real flag variety $M^\tau$
are described by the inequalities $w\ti\lambda\in\lambda-\ca C$, where
$w$ ranges over $W$.  By Lemma \ref{lemma;cone} this is equivalent to
$w\ti\lambda\in\lambda-\ca C^\a$ for all $w\in W$.  This implies
Kostant's result $\Delta(M^\tau)=\hull W^\a\lambda$, which was
discussed in Example \ref{example;kostant}.

\subsection*{Products of orbits}

Let us apply Theorem \ref{theorem;orbit} to the diagonal inclusion of
a group into two copies of itself.  We denote the ``small'' group by
$U$ and the ``large'' group by $U^2$.  We take the involution on $U^2$
to be the direct product of $\sigma$ with itself; the diagonal $f(U)$
is obviously stable under this involution.  The projection
$f^*\colon(\t^*)^2\to\t^*$ is the addition map and the homomorphism
$\phi^*\colon H^\bu(X,\Z)^{\otimes2}\cong H^\bu(X^2,\Z)\to
H^\bu(X,\Z)$ is the cup product.  Furthermore, the small Weyl chamber
$\t_+$ is contained diagonally in the big chamber $\t_+^2$ and hence
$W\rel=\{1\}$.  We obtain the following result from Theorem
\ref{theorem;orbit}.

\begin{theorem}\label{theorem;klyachko}
Let $\lambda$\upn, $\mu$ and $\nu\in\a^*_+$.  Then
$\nu\in\Delta(K\lambda\times K\mu)$ if and only if
\begin{equation}\label{equation;add}
w\inv\nu\le u\inv\lambda+v\inv\mu
\end{equation}
for all $u$\upn, $v$ and $w\in W$ such that $\sigma_w$ is contained in
$\sigma_u\cup\sigma_v$.
\qed
\end{theorem}

Here $\le$ denotes the partial ordering on $\a^*$ defined by the cone
$\ca C^\a$.  There is also a more ``economical'' version of this
result, in which the vector inequalities \eqref{equation;add} are
replaced by scalar inequalities, and the conditions on Schubert
classes for the general flag variety are replaced with conditions on
Schubert classes for the Grassmannians of $U$; cf. \cite[Theorem
4.2.1]{berenstein-sjamaar;coadjoint}.  See Example \ref{example;cross}
for a discussion of the $G$-orbits on $K\lambda\times K\mu$.  The Ky
Fan inequalities of Tam \cite{tam;unified-ky-fan} are a special case
of \eqref{equation;add}.

For $U=\SU(n)$, $\sigma=$ complex conjugation and $K=\SO(n)$ Theorem
\ref{theorem;klyachko} was also observed by Fulton as a corollary of a
result of Klyachko \cite{klyachko;stable-bundles-hermitian}.  In this
case we have $\a=\t$ and so the moment polytope of $M$ is the same as
that of $M^\tau$.  Theorem \ref{theorem;klyachko} then represents a
complete list of inequalities satisfied by the eigenvalues of a sum
$A+B$ of symmetric matrices in terms of the eigenvalues of $A$ and
$B$.  Such inequalities were first derived by Weyl.  See Fulton's
paper \cite{fulton;eigenvalues-sums} for a survey.

\subsection*{Isotropic flags}

We discuss Theorem \ref{theorem;klyachko} in more detail for
$U=\SU(n)$, where we take the involution to be the inner automorphism
defined by
$$
I_{pq}=
\begin{pmatrix}
-I_p&0\\0&I_q
\end{pmatrix}
$$
with $p+q=n$ and $p\le q$.  We start by giving an explicit description
of the varieties $K\lambda$.  Observe that $U^\sigma=K=\mathbf
S\bigl(\U(p)\times\U(q)\bigr)$, $G=\SU(p,q)$, and that $\q$ consists
of all matrices of the form
\begin{equation}\label{equation;rectangular}
\begin{pmatrix}
0&Z\\-Z^*&0
\end{pmatrix},
\end{equation}
where $Z$ is a complex $p\times q$-matrix.  The action of $K$ is given
by
\begin{equation}\label{equation;action}
(C,D)Z=CZD^*
\end{equation}
for $C\in\U(p)$ and $D\in\U(q)$.  The restricted root system is
$\BC_p$ if $p<q$ and $\CC_p$ if $p=q$.  (See e.g.\ \cite[Ch.\ VII,
\S2]{loos;symmetric-spaces;;1969}.)  Consider the antidiagonal
$p\times p$-matrices
\begin{equation}\label{equation;singular-values}
J=
\begin{pmatrix}
&&&1\\&&\antiddots&\\&1&&\\1&&&
\end{pmatrix},
\qquad
\Lambda=
\begin{pmatrix}
&&&\lambda_p\\&&\antiddots&\\&\lambda_2&&\\ \lambda_1&&&
\end{pmatrix}
\end{equation}
with $\lambda_1$, $\lambda_2,\dots$, $\lambda_p$ real.  Then a maximal
abelian subspace $\a$ of $\q$ is given by all $n\times n$-matrices of
the form
\vskip.5ex
$$
\lambda=
\begin{pmatrix}
\smash{\overset{p}{\vphantom{\Big(}0}}&
\smash{\overset{q-p}{\vphantom{\Big(}0}}& 
\smash{\overset{p}{\vphantom{\Big(}-J\Lambda J}}\\
0&0&0\\
\Lambda&0&0
\end{pmatrix}
\>
\begin{matrix}
{\scriptstyle p}\hfill\mbox{}\\{\scriptstyle q-p}\\{\scriptstyle
p}\hfill\mbox{}
\end{matrix},
$$
and $\a^*_+\cong\a_+$ is given by $\lambda_1\ge\lambda_2\ge\dots
\ge\lambda_p\ge0$.  The unitary matrix
$$
E=\frac1{\sqrt2}
\begin{pmatrix}
I&0&iJ\\0&\sqrt2\,I&0\\J&0&-iI
\end{pmatrix}
\in\U(n)
$$
diagonalizes $\a$ and maps $\a_+$ into the standard chamber $\t_+$,
$$
\Ad(E)\lambda=
\begin{pmatrix}
iJ\Lambda&0&0\\0&0&0\\0&0&-i\Lambda J
\end{pmatrix}.
$$
Now fix $\lambda\in\a^*_+$ and let $d=(d_1,d_2,\dots,d_k)$ be the
partition of $p$ determined by the conditions
$$
\lambda_1=\dots=\lambda_{d_1}>\lambda_{d_1+1}=\dots
=\lambda_{d_1+d_2}>\dots>\lambda_{d_1+\dots+d_{k-1}}=\dots=\lambda_p.
$$
Let $\ca F_d$ be the manifold of complex flags
$$
F\colon V_0\subset V_1\subset\cdots\subset V_{2k+1}\subset V_{2k+2}
$$
in $\C^n$ such that $V_0=\{0\}$, $V_{2k+2}=\C^n$ and
$$
\dim V_i-\dim V_{i-1}=
\begin{cases}
d_i&\text{if $1\le i\le k$}\\
q-p&\text{if $i=k+1$}\\
d_{2k+2-i}&\text{if $k+2\le i\le2k+1$.}
\end{cases}
$$
Then $\Ad(E)\lambda$ has the same stabilizer as the ``standard'' flag
$$
F_0\colon\{0\}\subset\C^{d_1}\subset\C^{d_1+d_2}\subset\cdots\subset
\C^p\subset\C^q\subset\C^{q+d_k}\subset\cdots\subset\C^{p+q-d_1}
\subset\C^n.
$$
It follows that there is a unique $U^\C$-equivariant isomorphism
$\psi\colon\ca F_d\to U\lambda$ which sends $E\inv F_0$ to $\lambda$.
The flag \emph{perpendicular} to $F\in\ca F_d$ is defined by
$$
F^\perp\colon\{0\}\subset V_{2k+1}^\perp\subset
V_{2k}^\perp\subset\cdots\subset V_1^\perp\subset\C^n,
$$
where the orthogonal complements are taken relative to the
nondegenerate Hermitian form $\langle v,w\rangle_{pq} =\langle
EI_{pq}E\inv v,w\rangle$, which has signature $q-p$.  Then
$F^\perp\in\ca F_d$, so $\perp$ defines an involution on $\ca F_d$ and
we obtain the diagram
\begin{equation}\label{equation;flag}
\begin{CD}
\ca F_d@>\psi>>U\lambda\\
@V{\perp}VV@VV{\tau}V\\
\ca F_d@>\psi>>U\lambda.
\end{CD}
\end{equation}

\begin{lemma}\label{lemma;isotropic}
The diagram \eqref{equation;flag} commutes.  Hence $K\lambda$ is
$G$-equivariantly diffeomorphic to the manifold of flags $\{0\}\subset
V_1\subset\cdots\subset V_k$ in $\C^n$ such that $V_i$ is
$\langle\cdot,\cdot\rangle_{pq}$-isotropic and $\dim V_i=
\sum_{j=1}^id_j$ for $i\le k$.
\end{lemma}

\begin{proof}
The first statement follows easily from the fact that $F_0=F_0^\perp$
and 
$$
\langle EI_{pq}E\inv Av,AEI_{pq}E\inv w\rangle_{pq} =\langle
v,w\rangle_{pq}
$$
for all $A\in\U(n)$.  To prove the second statement, recall that
$K\lambda=(U\lambda)^\tau$, which implies that $K\lambda$ is
isomorphic to the variety of self-perpendicular flags in $\ca F_d$.
Now observe that if $F=F^\perp$, then the truncated flag $F_{\le
k}\colon \{0\}\subset V_1\subset\cdots\subset V_k$ is isotropic and
$F$ is completely determined by $F_{\le k}$.
\end{proof}

Every matrix of the form \eqref{equation;rectangular} is in the
$K$-orbit of a unique matrix $\lambda$ in $\a^*_+$.  The numbers
$\lambda_1\ge\lambda_2\ge\dots \ge\lambda_p\ge0$ in
\eqref{equation;singular-values} are known as the \emph{singular
values}.  Let us call the vector
$(\lambda_1,\lambda_2,\dots,\lambda_p)$ corresponding to a $p\times
q$-matrix $Z$ the \emph{singular spectrum} of $Z$.  For $\lambda$ and
$\mu\in\a^*_+$ the polytope $\Delta(K\lambda\times K\mu)$ is the
collection of singular spectra of all complex matrices of the form
$A+B$, where $A$ and $B$ are $p\times q$-matrices with singular
spectra $\lambda$ and $\mu$, respectively.  (One gets a corresponding
result for singular values of \emph{real} $p\times q$-matrices by
considering $U=\SO(p+q)$ with the inner involution given by $I_{pq}$.
Although the abelian subspace $\a$ is the same here as in the complex
case, the root systems are different and one finds different
polytopes.)  We wish to thank W. Fulton for pointing out the
relationship with singular values.  See also his preprint
\cite{fulton;eigenvalues-invariant-factors}, which appeared after the
first version of this paper.

In the Lorentzian case $p=1$ it is easy to find the inequalities for
$\Delta(K\lambda\times K\mu)$ directly.  It is plain from
\eqref{equation;action} or Lemma \ref{lemma;isotropic} that the orbits
of $K$ on $\q$ are $2n-3$-spheres.  The restricted Weyl chamber is a
half-line, and the quotient map $\q^*\to\R_{\ge0}$ assigns to a sphere
its radius.  The moment ``polytope'' of a product of two spheres
$K\lambda$ and $K\mu$ is therefore simply the interval between
$\lvert\lambda-\mu\rvert$ and $\lambda+\mu$.  According to Example
\ref{example;cross}, there are two $G$-orbits on $S^{2n-3}\times
S^{2n-3}$ corresponding to the two elements of $W^\a$, namely the
diagonal and its complement.  The subpolytopes associated with these
orbits are $\{\lambda+\mu\}$, resp.\ the whole interval.

Let us also write out the case $p=q=2$.  The inequalities
\eqref{equation;add} can be enumerated explicitly as explained in
\cite{berenstein-sjamaar;coadjoint} or \cite{fulton;eigenvalues-sums}.
They are of the form
$$
\sum_{k\in\ca K}\nu_k\le\sum_{i\in\ca I}\lambda_i+\sum_{j\in\ca
J}\mu_j,
$$
where $\ca I\ca J\ca K$ range over all triples of ordered subsets of
$\{1,2,3,4\}$ listed in Table \ref{table;triple}.  For instance, the
triple $(1,2)(2,4)(2,4)$ encodes the inequality
$\nu_2+\nu_4\le\lambda_1+\lambda_2+\mu_2+\mu_4$.  For each triple $\ca
I\ca J\ca K$ there is a similar inequality corresponding to $\ca J\ca
I\ca K$.  To save space, for each $\ca I$, $\ca J$, $\ca K$ only one
of the triples $\ca I\ca J\ca K$ and $\ca J\ca I\ca K$ is listed in
the table.  Next to each triple $\ca I\ca J\ca K$ is listed the
\emph{dual} triple $\ca I^*\ca J^*\ca K^*$.  (If $\ca
I=(i_1<i_2<\dots<i_k)$ is a $k$-tuple of integers between $1$ and $n$,
and $\{1,2,\dots,n\}\setminus\ca I =(i_{k+1}<i_{k+2}<\dots<i_n)$ is
the complementary $n-k$-tuple, then $\ca I^*$ is defined as
$(n+1-i_n,n+1-i_{n-1},\dots,n+1-i_{k+1})$.  The inequality associated
with $\ca I^*\ca J^*\ca K^*$ is obtained from the one associated with
$\ca I\ca J\ca K$ by substituting $\lambda\to\lambda^*$,
$\mu\to\mu^*$, $\nu\to\nu^*$, where $\lambda^*=-w_0\lambda$.)  For
those triples which are self-dual we have omitted the corresponding
entry in the second column.

\begin{table}
\caption{Horn-Klyachko inequalities for $\SU(4)$}
\label{table;triple}
\begin{center}
\begin{tabular}[t]{|c|c|}
\hline
\emph{triple}&\emph{dual triple}\\
\hline\hline
(1)(1)(1)&(1,2,3)(1,2,3)(1,2,3)\\
(1)(2)(2)&(1,2,3)(1,2,4)(1,2,4)\\
(1)(3)(3)&(1,2,3)(1,3,4)(1,3,4)\\
(1)(4)(4)&(1,2,3)(2,3,4)(2,3,4)\\
(2)(2)(3)&(1,2,4)(1,2,4)(1,3,4)\\
(2)(3)(4)&(1,2,4)(1,3,4)(2,3,4)\\
\hline
\end{tabular}
\quad
\begin{tabular}[t]{|c|c|}
\hline
\emph{triple}&\emph{dual triple}\\
\hline\hline
(1,2)(1,2)(1,2)&\\
(1,2)(1,3)(1,3)&\\
(1,2)(1,4)(1,4)&(1,2)(2,3)(2,3)\\
(1,2)(2,4)(2,4)&\\
(1,2)(3,4)(3,4)&\\
(1,3)(1,3)(1,4)&(1,3)(1,3)(2,3)\\
(1,3)(1,4)(2,4)&(1,3)(2,3)(2,4)\\
(1,3)(2,4)(3,4)&\\
(1,4)(1,4)(3,4)&(2,3)(2,3)(3,4)\\
\hline
\end{tabular}
\end{center}
\end{table}

\begin{figure}
\setlength{\unitlength}{0.15mm}
$$
\begin{picture}(430,300)(0,-60)
\thinlines
\dashline{4.000}(0,0)(420,0)
\dashline{4.000}(0,0)(200,200)
\path(0,0)(0,60)
\path(2,52)(0,60)(-2,52)
\path(0,0)(60,-60)
\path(52.929,-55.757)(60,-60)(55.757,-52.929)
\texture{c0c0c0c0 0 0 0 0 0 0 0 c0c0c0c0 0 0 0 0 0 0 0 c0c0c0c0 0 0 0
	0 0 0 0 c0c0c0c0 0 0 0 0 0 0 0 }
\shade\path(300,90)(330,60)(390,60)
	(420,90)(420,150)(390,180)
	(330,180)(300,150)(300,90)
\texture{cccccccc 0 0 0 cccccccc 0 0 0 cccccccc 0 0 0 cccccccc 0 0 0
	cccccccc 0 0 0 cccccccc 0 0 0 cccccccc 0 0 0 cccccccc 0 0 0 }
\shade\path(390,180)(390,60)(420,90)(420,150)(390,180)
\thicklines
\path(390,180)(390,60)(420,90)(420,150)(390,180)
\put(0,0){\circle*{5}}
\put(60,0){\circle*{5}}
\put(30,30){\circle*{5}}
\put(60,30){\circle*{5}}
\put(360,120){\circle*{5}}
\put(420,150){\circle*{5}}
\put(60,-20){\makebox(0,0)[lb]{\smash{$\pi_1$}}}
\put(20,50){\makebox(0,0)[lb]{\smash{$\pi_2$}}}
\put(75,-60){\makebox(0,0)[lb]{\smash{$\alpha_1$}}}
\put(0,75){\makebox(0,0)[lb]{\smash{$\alpha_2$}}}
\put(75,30){\makebox(0,0)[lb]{\smash{$\lambda$}}}
\put(345,105){\makebox(0,0)[lb]{\smash{$\mu$}}}
\put(435,150){\makebox(0,0)[lb]{\smash{$\lambda+\mu$}}}
\end{picture}
$$
\caption{Moment polytopes of $K\lambda\times K\mu$ and $\bar\ca O_w$
with $K=\mathbf S\bigl(\U(2)\times\U(2)\bigr)$,
$\lambda=\frac1{2}\pi_1+\pi_2$, $\mu=4(\pi_1+\pi_2)$, $w=s_2s_1$}
\label{figure;su22}
\end{figure}
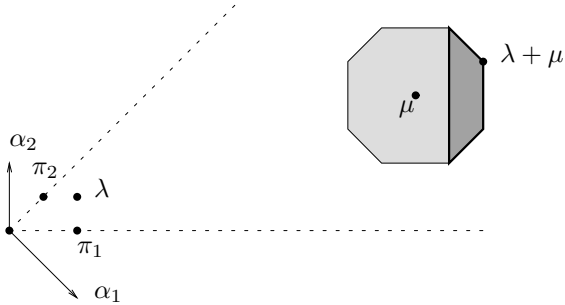

We get the inequalities for $\Delta(M^\tau)$ by substituting
$\lambda_1+\lambda_4=0$ and $\lambda_2+\lambda_3=0$ (plus similar
equalities for $\mu$ and $\nu$) into the inequalities for $\Delta(M)$.
Most inequalities turn out to be redundant.  In particular, every
inequality becomes \emph{identical} to its dual inequality, and all
non-self-dual inequalities with $\#\ca I=\#\ca J=\#\ca K=2$ drop out.
Here is a complete list, consisting of eighteen inequalities (besides
the obvious inequalities $\lambda_1\ge\lambda_2\ge0$,
$\mu_1\ge\mu_2\ge0$, $\nu_1\ge\nu_2\ge0$):
\begin{xalignat*}{2}
\nu_1 &\le\lambda_1+\mu_1 
&\nu_1+\nu_2 &\le\lambda_1+\lambda_2+\mu_1+\mu_2\\
\nu_2 &\le
\begin{cases}
\lambda_1+\mu_2\\ \lambda_2+\mu_1
\end{cases}
&\nu_1-\nu_2 &\le
\begin{cases}
\lambda_1+\lambda_2+\mu_1-\mu_2\\ \lambda_1-\lambda_2+\mu_1+\mu_2
\end{cases}\\
\nu_2 &\ge
\begin{cases}
-\lambda_1+\mu_2\\ \lambda_2-\mu_1
\end{cases}
&\nu_1-\nu_2 &\ge
\begin{cases}
-\lambda_1-\lambda_2+\mu_1-\mu_2\\ \lambda_1-\lambda_2-\mu_1-\mu_2
\end{cases}\\
\nu_1 &\ge\lvert\lambda_1-\mu_1\rvert
&\nu_1+\nu_2 &\ge\lvert\lambda_1+\lambda_2-\mu_1-\mu_2\rvert\\
\nu_1 &\ge\lvert\lambda_2-\mu_2\rvert
&\nu_1+\nu_2 &\ge\lvert\lambda_1-\lambda_2-\mu_1+\mu_2\rvert.
\end{xalignat*}
See Figure \ref{figure;su22} for an example.  Also shown, in dark, is
the subpolytope $\Delta(\bar\ca O_w)$, where $w$ is the product of the
simple reflections $s_1$ and $s_2$ in $W^\a$.  The simple roots are
denoted by $\alpha_1$ and $\alpha_2$ and the corresponding fundamental
weights by $\pi_1$ and $\pi_2$.

\subsection*{An unsolved problem}

As is obvious from the above examples, the inequalities of Theorem
\ref{theorem;orbit} are in general highly overdetermined.  We believe
it should be possible to find a complete set of inequalities which
uses only data involving the restricted Weyl groups and cohomology
classes on real flag varieties.

\begin{problem}\label{problem;real}
Does the statement of Theorem {\rm\ref{theorem;orbit}} remain true if
we replace $\ti W$\upn, $W$\upn, and $W\rel$ by $\ti W^{\ti\a}$\upn,
$W^\a$\upn, and $W\rel^\a$\upn, respectively\upn, and complex Schubert
classes by $\Z/2$-Schubert classes on appropriate real flag
varieties\upn?
\end{problem}

Here $\ti W^{\ti\a}$ and $W^\a$ are the restricted Weyl groups.  The
restricted relative Weyl set $W\rel^\a$ is defined in the same way as
$W\rel$, except that in the definition $W$ is replaced by $W^\a$,
$\ti\t$ by $\ti\a$, and $\t$ by $\a$.

In the example of the maximal torus the answer to Problem
\ref{problem;real} is affirmative.  In the $K\lambda\times
K\mu$-example it is not very hard to see that the same is true in the
low-dimensional cases which we wrote out, although in general the
question appears to be nontrivial.

\appendix

\section{Anti-symplectic maps}\label{section;antisymplectic}

The results in this appendix can be found in a somewhat weaker form in
the papers of Meyer \cite{meyer;hamiltonian-discrete} and Duistermaat
\cite{duistermaat;convexity}.  Let $(V,\omega)$ be a symplectic vector
space.  Recall that a complex structure $J$ on $V$ is
\emph{compatible} with $\omega$ if $J\in\Sp(V,\omega)$ and the
symmetric bilinear form $\omega(\cdot,J\cdot)$ is positive definite.
Let $\ca J(\omega)$ be the set of complex structures on $V$ that are
compatible with $\omega$ and $\Inn(V)$ the set of (positive definite)
inner products on $V$.  Recall that we have maps $\iota\colon\ca
J(\omega)\to\Inn(V)$ and $\pi\colon\Inn(V)\to\ca J(\omega)$ defined by
$$
\iota(J)(v,w)=\omega(v,Jw),\qquad
\pi(\beta)=(-A_\beta^2)^{-1/2}A_\beta,
$$
where for $\beta\in\Inn(V)$ the antisymmetric linear map $A_\beta$ is
defined by $\omega(v,w)=\beta(A_\beta v,w)$.  Furthermore, the map
$\iota$ is injective and $\pi$ is a left-inverse of $\iota$.  It
follows that $\pi$ is a deformation retraction and that $\ca
J(\omega)$ is contractible.  (See McDuff and Salamon
\cite[\S2.5]{mcduff-salamon;introduction}.)

Now let $\ca G$ be a Lie group with a real character $\eps\colon\ca
G\to\R^\times$ and suppose that $\ca G$ acts linearly and
\emph{sesquisymplectically} on $V$ in the sense that
$g^*\omega=\eps(g)\omega$ for $g\in\ca G$.  Define an action of $\ca
G$ on $\ca J(\omega)$ by
$$
g\cdot J=\eps(g)g\circ J\circ g\inv.
$$
Then $J$ is invariant for this action if $J\circ g=\eps(g)g\circ J$.
Let $\ca G$ act on $\Inn(V)$ by $(g\cdot\beta)(v,w)=\beta(gv,gw)$.
The proof of the following result is straightforward.

\begin{lemma}\label{lemma;complex}
The maps $\iota$ and $\pi$ are $\ca G$-equivariant.  Hence\upn, if
$\ca G$ is compact\upn, $\ca J(\omega)^{\ca G}$ is contractible and
nonempty.
\qed
\end{lemma}

Of course, if $\ca G$ is compact, $\eps$ takes values in $\{\pm1\}$.
Let us denote the kernel of $\eps$ by $\ca G_1$.

\begin{lemma}\label{lemma;realeigen}
Assume that $\ca G$ is compact.  Let $g\in\ca G$ be antisymplectic.
Then the following conditions are equivalent.
\begin{enumerate}
\item\label{item;ginvolution}
$g$ is an involution.
\item\label{item;realeigen}
All eigenvalues of $g$ are real.
\item\label{item;subspace}
$V^g$ is a Lagrangian subspace.
\end{enumerate}
If any of these conditions is satisfied\upn, then $\ca G\cong\ca
G_1\rtimes\{\pm1\}$ and $V^g$ is totally real with respect to any $\ca
G$-invariant almost complex structure.
\end{lemma}

\begin{proof}
It is obvious that \eqref{item;ginvolution} implies
\eqref{item;realeigen}.

Assume that \eqref{item;realeigen} holds.  According to Lemma
\ref{lemma;complex} there exists $J\in\ca J(\omega)^{\ca G}$.  Then
$gJ=-Jg$.  Moreover, the eigenvalues of $g$ are $\pm1$ because $g$ is
orthogonal relative to the inner product $\iota(J)$.  Let
$\lambda=\pm1$ be an eigenvalue of $g$ with eigenvector $e\in V$ and
put $f=Je$.  Then
\begin{equation}\label{equation;pm}
gf=gJe=-Jge=-\lambda f=\mp f.
\end{equation}
Let $W$ be the real span of $e$ and $f$.  Then $W$ is a Hermitian
subspace of $W$.  Since $g$ is antisymplectic and preserves $W$, it
preserves the symplectic subspace $W^\omega=W^\perp$.  By induction on
the dimension we may assume that the fixed points of $g|_{W^\perp}$
form a Lagrangian subspace.  On the other hand, it is clear that
$g|_W$ fixes either $e$ or $f$, depending on the value of $\lambda$.
This implies \eqref{item;subspace}.

Now suppose that \eqref{item;subspace} holds.  Let
$\{e_1,e_2,\dots,e_n\}$ be a basis of $V^g$ and put $f_i=Je_i$.  It
follows from \eqref{equation;pm} that $f_i$ has eigenvalue $-1$ and
therefore
$$
\{e_1,f_1,e_2,f_2,\dots,e_n,f_n\}
$$
is a (symplectic) basis of $V$.  Hence \eqref{item;ginvolution} holds.

If $g\in\ca G\setminus\ca G_1$ is an involution, the exact sequence
$\ca G_1\hookrightarrow\ca G\to\kern-1.9ex\to\{\pm1\}$ splits and
hence $\ca G\cong\ca G_1\rtimes\{\pm1\}$.  The fact that $V^g$ is
totally real follows from the following trivial observation: let $W$
be a (real) vector space with a complex structure $J$ and let $\tau$
be an involution which is antiholomorphic in the sense that $\tau
J=-J\tau$.  Then the subspaces $W^\tau$ and $W^{-\tau}=JW^\tau$ are
totally real, $\tau$ is complex conjugation with respect to $W^\tau$,
and $W=W^\tau\oplus JW^\tau\cong(W^\tau)^\C$.
\end{proof}

Now let $(M,\omega)$ be a connected symplectic manifold and assume
that $\ca G$ acts sesquisymplectically on $M$, i.e.\
$g^*\omega=\eps(g)\omega$ for $g\in\ca G$.  Let us denote by $\ca
J(M,\omega)$ the space of almost complex structures on $M$ which are
compatible with $\omega$.  As before, $\ca G$ acts in a natural way on
$\ca J(M,\omega)$.

\begin{lemma}\label{lemma;almostcomplex}
Assume that $\ca G$ acts properly on $M$.
\begin{enumerate}
\item\label{item;complex}
$\ca J(M,\omega)^{\ca G}$ is nonempty and connected.
\item\label{item;lagrangian}
If there exist $g\in\ca G\setminus\ca G_1$ and $m\in M$ such that $g$
fixes $m$ and has only real eigenvalues at $m$\upn, then $g$ is an
involution and $M^g$ is a Lagrangian submanifold.  Moreover\upn, $M^g$
is totally real with respect to any $\ca G$-invariant almost complex
structure.
\end{enumerate}
\end{lemma}
  
\begin{proof}
\eqref{item;complex} follows from Lemma \ref{lemma;complex} by a
standard argument as in \cite[\S2.6]{mcduff-salamon;introduction}.
\eqref{item;lagrangian} is proved by applying Lemma
\ref{lemma;realeigen} to the stabilizer subgroup $\ca G_m$ (which is
compact and contains $g$) acting on the tangent space $T_mM$.
\end{proof}

\section{Symmetric pairs}\label{section;symmetric}

This appendix is a collection of elementary facts concerning the
symmetric pair $(U,K)$.  In large part, this is well-known material
culled from Helgason \cite{helgason;differential-geometry-lie-groups}
and Loose \cite{loos;symmetric-spaces;;1969}.  Let $U$ be a compact
connected Lie group with involution $\sigma$ and let $K=(U^\sigma)_0$.
Let
$$
Q=\{\,u\sigma(u)\inv\mid u\in U\,\}
$$
be the submanifold of symmetric elements.  Put $\q=T_1Q$.  Then
$Q=\exp\q$ and the map $U\to Q$ which sends $u$ to $u\sigma(u)\inv$
induces an isometry from the symmetric space $U/U^\sigma$ onto $Q$.
(See e.g.\ Chapters II and VI in \cite{loos;symmetric-spaces;;1969}.)
The Lie algebra of $U$ decomposes into the $\pm1$-eigenspaces for the
involution $\sigma$,
\begin{equation}\label{equation;decompose}
\u=\k\oplus\q
\end{equation}
and dually we have $\u^*=\k^*\oplus\q^*$.  Let $A$ be a maximal torus
(connected flat submanifold) in $Q$ containing $1$ and let $\a=T_1A$.
Then $\a$ is a maximal abelian subspace of $\q$ and $A=\exp\a$.
Clearly $A$ is stable under $\sigma$.  The fixed-point group (``real
part'') is
$$
A^\sigma=A[2]=\{\,a\mid a^2=1\,\},
$$
the subgroup of $2$-torsion elements.  Hence $\#A^\sigma=2^r$, where
$r=\dim A$, the rank of the symmetric space $Q$.  The following
assertion is a straightforward corollary of the $KAK$-decomposition,
Theorem 8.6 in \cite[Ch.\
VII]{helgason;differential-geometry-lie-groups}.

\begin{lemma}\label{lemma;fixed}
$U^\sigma=KA[2]=A[2]K$.
\qed
\end{lemma}

The factors $K$ and $A[2]$ usually have a nontrivial intersection, so
this product is rarely semidirect.  Now choose a maximal torus $T$ of
$U$ such that $T\supset A$.  Such a torus is always $\sigma$-stable.
Let $K'=\ca Z_K(A)$, the centralizer of $A$ in $K$.  Then
$T'=(K')_0\cap T$ is a maximal torus of (the identity component of)
$K'$.  Moreover,
\begin{equation}\label{equation;torus}
T=T'A
\end{equation}
and $T'\cap A\subset A[2]$.  (See Proposition 3.3 in \cite[Ch.\
VI]{loos;symmetric-spaces;;1969}.)  From \eqref{equation;decompose}
and \eqref{equation;torus} we get the direct-sum decomposition
\begin{equation}\label{equation;cartandecompose}
\t=\t'\oplus\a.
\end{equation}
Let $R$ be the root system of $(U,T)$.  Since $\sigma$ is an
automorphism and $\t$ is $\sigma$-stable, $\sigma$ acts on $R$.
Define
$$
R'=\{\,\alpha\in R\mid\sigma(\alpha)=\alpha\,\}, \qquad
R^\a=\{\,\alpha|_\a\mid\alpha\not\in R'\,\}.
$$
Then $R'$ is the set of roots which vanish on $\a$.  It is a root
subsystem of $R$, namely the root system of $(K',T')$, and $R^\a$ is a
(not necessarily reduced) root system in $\a^*$, known as the system
of \emph{restricted roots} of the symmetric pair $(\u,\k)$.  We define
sets of positive roots in $R$ and $R^\a$ as follows: fix sets of
positive roots $R'_+$ in $R'$ and $R^\a_+$ in $R^\a$.  Then
\begin{equation}\label{equation;positive}
R_+=R'_+\cup\{\,\alpha\in R\mid\alpha|_\a\in R^\a_+\,\}
\end{equation}
is a system of positive roots in $R$.  Let
$$
\t_+=\{\,\xi\in\t\mid\alpha(\xi)\ge0\;\text{for $\alpha\in R_+$}\,\}
\quad\text{and}\quad \a_+=\{\,\xi\in\t\mid\alpha(\xi)\ge0\;\text{for
$\alpha\in R^\a_+$}\,\}
$$
be the associated Weyl chambers in $\t$, resp.\ $\a$.  Then $\t_+$ is
a fundamental domain for the $U$-action on $\u$ and $\a_+$ is a
fundamental domain for the $K$-action on $\q$.  Let
$\t^*_+\subset\t^*$ and $\a^*_+\subset\a^*$ be the corresponding dual
Weyl chambers.  Because of the way the positive roots in $\t$ and $\a$
are related we have the following equalities.

\begin{lemma}\label{lemma;chambers}
$\a_+=\t_+\cap\q$ and $\a^*_+=\t^*_+\cap\q^*$.
\qed
\end{lemma}

Likewise, let $\ca C$ be the \emph{root cone}, i.e.\ the cone in
$\t^*$ spanned by the positive roots $R_+$, and let $\ca C^\a$ be the
\emph{restricted root cone}, i.e.\ the cone in $\a^*$ spanned by
$R^\a_+$.  Then $\ca C$ is the dual cone of $\t_+$ and $\ca C^\a$ is
the dual cone of $\a_+$.  The root cones are related as follows.

\begin{lemma}\label{lemma;cone}
$\ca C^\a=\ca C\cap\a^*$.
\end{lemma}

\begin{proof}
The inclusion $\ca C^\a\subset\ca C\cap\a^*$ is obvious.  The reverse
inclusion follows easily from the fact that
$\alpha(\xi)=\frac1{2}\bigl(\alpha+\sigma(\alpha)\bigr)(\xi)$ for
$\alpha\in R$ and $\xi\in\a$.
\end{proof}

Finally we recall the relationship between the Weyl groups of the root
systems $R$, $R'$ and $R^\a$, which we denote by $W$, $W'$ and $W^\a$,
respectively.  Observe that $\sigma$ acts in a natural way on $W$
because $\t$ is $\sigma$-stable.

\begin{lemma}\label{lemma;weyl}
\begin{enumerate}
\item\label{item;normal}
$W'\cong\ca Z_W(\a)$\upn, $W^\sigma=\ca N_W(\a)$\upn, and $W^\a\cong
W^\sigma/W'$.
\item\label{item;suborbit}
$(U\lambda)^A\cap\q^*=U\lambda\cap\a^* =W\lambda\cap\a^*=W^\a\lambda$
for all $\lambda\in\a^*$.
\item\label{item;hull}
$(\hull W\lambda)\cap\a^*=\hull W^\a\lambda$ for all $\lambda\in\a^*$.
\end{enumerate}
\end{lemma}

Here ``$\ca Z$'' denotes the centralizer of a subalgebra, ``$\ca N$''
its normalizer, and ``$\hull$'' denotes the convex hull of a subset of
a real vector space.

\begin{proof}
The isomorphisms $W'\cong\ca Z_W(\a)$ and $W^\a\cong\ca N_W(\a)/\ca
Z_W(\a)$ are in
\cite[p. 325]{helgason;differential-geometry-lie-groups}.  If $w\in
W^\sigma$, then $\sigma(w\xi)=\sigma(w)\sigma(\xi)=-w\xi$ for all
$\xi\in\a$, and therefore $w\in\ca N_W(\a)$.  Conversely, if $w\in\ca
N_W(\a)$, then $w$ preserves the decomposition
\eqref{equation;cartandecompose}.  (This follows for instance from the
fact that the decompositions \eqref{equation;decompose} and
\eqref{equation;cartandecompose} are orthogonal with respect to the
Killing form.)  Now let $\xi\in\t$ and write $\xi=\xi^++\xi^-$ with
$\xi^+\in\t'$ and $\xi^-\in\a$.  Then
$$
\sigma(w)\xi =\sigma(w)(\xi^++\xi^-)
=\sigma\bigl(w\sigma(\xi^+)+w\sigma(\xi^-)\bigr)
=\sigma(w\xi^+-w\xi^-) =w\xi^++w\xi^-=w\xi
$$
and hence $\sigma(w)=w$.  The conclusion is that $\ca
N_W(\a)=W^\sigma$.  This proves \eqref{item;normal}.

The first equality in \eqref{item;suborbit} follows easily from the
fact that $A$ is maximal abelian in $Q$.  To prove the second
equality, observe that $U\lambda\cap\a^*\subset
U\lambda\cap\t^*=W\lambda$, so $U\lambda\cap\a^*\subset
W\lambda\cap\a^*$.  The last equality is Corollary 8.9 of \cite[Ch.\
VII]{helgason;differential-geometry-lie-groups}.

The inclusion $\hull W^\a\lambda\subset(\hull W\lambda)\cap\a$ follows
from \eqref{item;suborbit}.  The reverse inclusion follows from the
fact that
$$
\hull W\lambda=\bigcap_{w\in W}w(\lambda-\ca C)\qquad
\text{and}\qquad\hull W^\a\lambda=\bigcap_{w\in W^\a}w(\lambda-\ca
C^\a),
$$
combined with Lemma \ref{lemma;cone}.
\end{proof}

\section{Notation}\label{section;notation}

\begin{tabbing}
\indent \= $M_\mu$; $M_0=M\qu G$ \= \kill 
\> $U$; $\sigma$ \> compact connected Lie group; involution\\
\> $U^\sigma$; $K$ \> fixed-point group; its identity component\\
\> $Q$; $\q$ \> space of symmetric elements $u\sigma(u)\inv$; its
tangent space at $1$\\
\> $\u=\k\oplus\q$ \> eigenspace decomposition of $\u$ under
$\sigma$\\
\> $G$; $\theta$ \> real reductive group dual to $U$; Cartan
involution\\
\> $\g=\k\oplus\p$ \> Cartan decomposition of $\g$, where $\p=i\q$\\
\> $A$; $\a$ \> maximal torus in $Q$; its tangent space at $1$\\
\> $A[2]$; $K'$ \> $2$-torsion elements in $A$; centralizer of $A$ in
$K$\\
\> $T'$; $T=T'A$ \> maximal torus of $K'$; maximal torus of $U$\\
\> $\t'$; $\t=\t'\oplus\a$ \> Cartan subalgebra of $\k'$; Cartan
subalgebra of $\u$\\
\> $R$; $R'$; $R^\a$ \> roots of $(\u,\t)$; roots of $(\k',\t')$;
restricted roots\\
\> $W$; $W'$; $W^\a$ \> Weyl groups of $R$; $R'$; $R^\a$\\
\> $w_0$; $w'_0$; $w_0^\a$ \> longest elements of $W$; $W'$; $W^\a$\\
\> $\t^*_+$; $\a^*_+$ \> fundamental Weyl chambers for $R$; resp.\
$R^\a$\\
\> $\ca C$; $\ca C^\a$ \> cone spanned by $R_+$; resp.\ $R^\a_+$\\
\> $\Lambda^*_+$; $\sigma_+$ \> set of dominant weights; involution on
$\t^*_+$ induced by $\sigma$\\
\> $V_\lambda$ \> irreducible $U$-module with highest weight
$\lambda\in\Lambda^*_+$\\
\> $\gamma_\u$; $\gamma_\g$ \> complex conjugations on $\u^\C$
relative to $\u$; resp.\ $\g$\\
\> $(M,\omega)$ \> symplectic manifold acted on by $U$\\
\> $\Phi$; $\tau$; $M^\tau$ \> moment map; antisymplectic involution;
fixed-point set\\
\> $\xi_M$; $\Phi^\xi$ \> vector field induced by $\xi\in\u$;
component of $\Phi$ along $\xi$\\
\> $\Delta(M)$; $\Delta(M^\tau)$ \> moment polytopes
$\Phi(M)\cap\t^*_+$; resp.\ $\Phi(M^\tau)\cap\t^*_+$
\end{tabbing}


\bibliographystyle{amsplain}

\providecommand{\bysame}{\leavevmode\hbox to3em{\hrulefill}\thinspace}


\end{document}